\documentclass[a4paper,11pt]{article}
\usepackage{amsmath,amstext}  
\usepackage{amsfonts,amssymb, bbm} 
\usepackage[LGR, T1]{fontenc}
\usepackage[utf8]{inputenc}
\usepackage[english]{babel}
\usepackage{color}
\usepackage{graphicx, graphics}
\usepackage{subfigure}
\usepackage{epsfig}
\usepackage{epsf}
\usepackage{vmargin}
\newcommand{\ind}[1]{\mathbbm{1}_{\{#1\}}}
\newcommand{\wh}{\widehat}

\newcommand{\wt}{\widetilde}
\newcommand{\dE}{\mathbb{E}}
\newcommand{\dP}{\mathbb{P}}
\newcommand{\dR}{\mathbb{R}}
\newtheorem{theorem}{Theorem}[section]
\newtheorem{lemma}[theorem]{Lemma}
\newtheorem{rem}[theorem]{Remark}

\newtheorem{hyp}[theorem]{Assumption}
\newtheorem{corollary}[theorem]{Corollary}
\title{Approximate Kalman--Bucy filter for continuous-time semi-Markov jump linear systems}
\author{Beno\^\i te~de~Saporta \and Eduardo~F.~Costa
\thanks{Beno\^{\i}te de Saporta is with Univ. Montpellier 2
F-34095 Montpellier, France, 
CNRS I3M UMR 5149 F-34095 Montpellier, France
and INRIA Bordeaux Sud Ouest, team CQFD, F-33400 Talence, France. 
e-mail: Benoite.de-Saporta@univ-montp2.fr}
\thanks{Eduardo F. Costa is with Univ. S\~{a}o Paulo -
Instituto de Ci\^{e}ncias Mathem\'{a}ticas e de Computa\c{c}\~{a}o, C.P. 668,
13560-970, S\~{a}o Carlos, SP, Brazil. email: efcosta@icmc.usp.br}}
\begin{document}

\maketitle

\begin{abstract}
The aim of this paper is to propose a new numerical approximation of the 
Kalman--Bucy filter for semi-Markov jump linear systems.
This approximation is based on the selection of typical trajectories 
of the driving semi-Markov chain of the process by 
using an optimal quantization technique. 
The main advantage of this approach is that it makes pre-computations possible.
We derive a Lipschitz property for the solution of the Riccati equation 
and a general result on the convergence of perturbed solutions 
of semi-Markov switching Riccati equations 
when the perturbation comes from the driving semi-Markov chain.
Based on these results, we prove the convergence of our approximation scheme
in a general infinite countable state space framework 
and derive an error bound in terms of the quantization error and time discretization step.
We employ the proposed filter in a magnetic levitation example 
with markovian failures and compare its performance with both 
the Kalman--Bucy filter and the Markovian linear minimum mean squares estimator.
\end{abstract}

\section*{IEEE Copyright Notice}
© IEEE. Personal use of this material is permitted. However, permission to reprint/republish this material for advertising or promotional purposes or for creating new collective works for resale or redistribution to servers or lists, or to reuse any copyrighted component of this work in other works must be obtained from the IEEE.

This material is presented to ensure timely dissemination of scholarly and technical work. Copyright and all rights therein are retained by authors or by other copyright holders. All persons copying this information are expected to adhere to the terms and constraints invoked by each author's copyright. In most cases, these works may not be reposted without the explicit permission of the copyright holder.

For more details, see the IEEE Copyright Policy\\
\footnotesize \verb+http://www.ieee.org/publications_standards/publications/rights/copyrightpolicy.html+
\normalsize

\section{Introduction}
{M}{arkov} jump linear systems (MJLS) have been largely studied and 
disseminated during the last decades. MJLS have a relatively simple structure 
that allows for useful, strong properties 
\cite{CostaFragosoMarques05,CostaFragosoTodorov,Dragan2009,Dragan2013},
and 
provide suitable models for applications \cite{doVal99JE,Sworder83,Siqueira04},
with a booming field in web/internet based control \cite{Geromel09,Huang13}.
One limitation of MJLS is that the sojourn times between jumps is a time-homogeneous 
exponential random variable, thus motivating the study of a wider class of systems with 
general sojourn-time distributions, the so-called semi-Markov jump linear 
systems (sMJLS) or sojourn-time-dependent MJLS \cite{Huang13,Campo91,Schwartz03,Hou06,Huang13IJRNC}.

In this paper, we consider continuous-time sMJLS with instantaneous 
(or close to instantaneous) observation of the state of the semi-Markov 
chain at time instant $t$, denoted here by $\theta(t)$.
The state space of the semi-Markov chain may be infinite. 
We seek for an approximate optimal filter for the variable $x(t)$ 
that composes the state of the sMJLS jointly with $\theta(t)$. 
Of course, estimating the state component $x(t)$ is highly relevant 
and allows the use of standard control strategies like linear state feedback.

It is well known that the optimal estimator for $x(t)$ is given by the 
standard Kalman--Bucy filter (KBF) \cite{Anderson79,Jazwinski70,Kalman60,Kalman61,Kumar86}
because, given the observation of the past values of $\theta$, the 
distribution of the random variable $x(t)$ is exactly the same as 
in a time varying system. 
The main problem faced when implementing the KBF for MJLS or sMJLS, particularly 
in continuous time, is the pre-computation.
Pre-computation refers to the computation of the relevant parameters of the KBF 
and storage in the controller/computer memory prior to the system operation, 
which makes the implementation of the filter 
fast enough to couple with a wide range of applications. 
Unfortunately, pre-computation is not viable for (s)MJLS in continuous time,
as it involves solving a Riccati differential equation that branches at every 
jump time $T_k$, 
and the jumps can occur at any time instant according to an exponential distribution, 
so that pre-computation would involve computation of an infinite number of branches.
Another way to explain this drawback of the KBF is to say that the KBF is not a
Markovian linear estimator because the gain at time $t$ does not depend only on $\theta(t)$ 
but on the whole trajectory $\{\theta(s),0\leq s\leq t\}$.
This drawback of the KBF has motivated the development of other 
filters for MJLS, and one of the most successful ones is the  Markovian linear minimum 
mean squares estimator (LMMSE) that has been derived in \cite{FC10}, 
whose parameters can be pre-computed, 
see also \cite{CostaFragosoTodorov,Costa11autom_filter}. 
To our best knowledge, there is no pre-computable filter for sMJLS.

The filter proposed here is built in several steps. The first step is the discretization 
by quantization of the Markov chain, providing a finite number of typical trajectories.  
The second step consists in solving the Riccati differential equation on each of these trajectories 
and store the results. To compute the filter in real time, one just needs to select the appropriate 
pre-computed branch at each jump time and follow it until the next jump time. 
This selection step is made by looking up the projection of the real jump time in the quantization 
grid and choosing the corresponding Riccati branch. 
In case the real jump time is observed with some delay (non-instantaneous observation of $\theta$), 
then the observed jump time is projected in the quantization grid instead, see Remarks \ref{rem-byproduct}, \ref{rem-delayed-observation}.

The quantization technique selects optimized typical trajectories of the semi-Markov chain.  
Optimal quantization methods have been developed recently in numerical probability, 
nonlinear filtering or optimal stochastic control for diffusion processes with applications in finance \cite{bally03,bally05,pages98,pages05,pages04b,pages04}  or for piecewise deterministic Markov processes with applications in reliability \cite{brandejsky12,brandejsky13,saporta10,saporta12}. To our best knowledge, this technique has not been applied to MJLS or sMJLS yet. The optimal quantization of a random variable $X$ consists in finding a finite grid such that the projection $\widehat{X}$ of $X$ on this grid minimizes some $L^{p}$ norm of the difference $X-\widehat{X}$. Roughly speaking, such a grid will have more points in the areas of high density of $X$. 
One interesting feature of this procedure is that the construction of the optimized grids using 
the CLVQ algorithm (competitive learning vector quantization) \cite{pages98,gray98} 
only requires a simulator of the process and no special knowledge about the distribution of $X$. 

As explained for instance in \cite{pages04}, for the convergence of the quantized 
process towards the original process, some Lipschitz-continuity conditions are needed,   
hence we start investigating the Lipschitz continuity of solutions of Riccati equations. 
Of course, this involves evaluating the difference of two Riccati solutions, 
which is not a positive semi-definite nor a negative-definite matrix, 
preventing us to directly use the positive invariance property of Riccati equations,  
thus introducing some complication in the analysis given in Theorem \ref{lem-1-Ric}.
A by product of our procedure is a general result on the convergence of perturbed 
solutions of semi-Markov switching Riccati equations, 
when the perturbation comes from the driving semi-Markov chain and can be either a random 
perturbation of the jump times or a deterministic delay, or both, see Remark \ref{rem-byproduct}.
Regarding the proposed filter, we obtain an error bound w.r.t. the exact KBF depending on the quantization error and time discretization step. 
It goes to zero when the number of points in the grids goes to infinity. 

The approximation results are illustrated and compared with the exact KBF and the 
LMMSE in the Markovian framework for a numerical example of a magnetic suspension system, 
confirming via Monte Carlo simulation 
that the proposed filter is effective for state estimation 
even when a comparatively low number of points in the discretization grids is considered.

The paper is organized as follows. Section \ref{sec-problem} presents the 
KBF and the sMJLS setup. The KBF approximation scheme is explained 
in Section \ref{sec-KBF}, and its convergence is studied in Section \ref{sec-convergence}.
The results are illustrated in a magnetic suspension system, see Section \ref{sec-example},
and some concluding remarks are presented in Section \ref{sec-conclusion}.

\section{Problem setting}
\label{sec-problem}
We start with some general notation. 
For $z,\hat z\in\dR$, $z\wedge \hat z=\min\{z,\hat z\}$ is the minimum between $z$ and $\hat z$.
For a vector $X=(x_1,\ldots,x_n)\in\dR^n$, $|X|$ denotes its Euclidean norm $|X|^2=\sum x_i^2$ and $X'$ denotes its transpose. 
Let $\mathcal{C}(n)$ be the set of $n\times n$ symmetric positive definite matrices and $I_{n}$ (or $I$ when there is no ambiguity) the identity matrix of size $n\times n$. 
For any two symmetric positive semi-definite matrices $M$ and $\wh{M}$, $M\geq \wh{M}$ means that $M-\wh{M}$ is 
positive semidefinite and $M >\wh{M}$ means that $M-\wh{M}\in\mathcal{C}(n)$.
Let $\lambda_{\min}(M)$ and $\lambda_{\max}(M)$ denote the lowest and highest eigenvalue of matrix $M\in\mathcal{C}(n)$ respectively.
For a matrix $M\in\dR^{n\times n}$, $M'$ is the transpose of $M$ and $\|M\|$ stands for its $L^2$ matrix norm $\|M\|^2=\lambda_{\max}(M'M)$.

Let $(\Omega, \mathcal{F}, \dP)$ be a probability space, $\dE$ denote the expectation with respect to $\dP$, and $Var(X)$ is the variance-covariance matrix of the random vector $X$. 
Let $\{\theta(t), t\geq 0\}$ be a semi-Markov jump process on the countable state space $\mathcal{S}$.
We denote by $F_i$ the cumulative distribution function of the sojourn time of $\theta$ in state $i$.
For a family $\{M_i, i\in\mathcal{S}\}$ of square matrices indexed by $\mathcal{S}$, we set $\|M\|_{\mathcal{S}}=\sup_{i\in\mathcal{S}}\|M_i\|\leq \infty$.

We consider a sMJLS satisfying 
\begin{equation*}
\left\{\begin{array}{rcl}
dx(t)&=&A_{\theta(t)}x(t)dt+E_{\theta(t)}dw(t),\\
dy(t)&=&C_{\theta(t)}x(t)dt+D_{\theta(t)}dv(t),\label{mjls}
\end{array}\right.
\end{equation*}
for $0\leq t\leq T$, where $T$ is a given time horizon, 
$\big(x(t),\theta(t)\big)\in\dR^{n_1}\times \mathcal S$ is the state process, 
$y(t)\in\dR^{n_2}$ is the measurement process, 
$\{w(t), 0\leq t\leq T\}$ and $\{v(t), 0\leq t\leq T\}$ are independent standard Wiener processes with respective 
dimensions $n_3$ and $n_4$, independent from $\{\theta(t), t\geq 0\}$, and 
$\{A_i, i\in\mathcal{S}\}$, $\{C_i, i\in\mathcal{S}\}$, $\{D_i, i\in\mathcal{S}\}$ and $\{E_i, i\in\mathcal{S}\}$
are families of matrices with respective size $n_1\times n_1$, $n_2\times n_1$, $n_2\times n_4$ and $n_1\times n_3$ such that $D_iD_i'>0$ is nonsingular for all $i$ (nonsingular measurement noise). 
 
We use two different sets of assumptions for the parameters of our problems. The first one is more restrictive but relevant for applications, and the second more general one will be used in the convergence proofs.
\begin{hyp}\label{hyp:finite}
The state space $\mathcal{S}$ is finite, $\mathcal{S}=\{1,2,\ldots,N\}$ and the cumulative distribution functions of the sojourn times $F_i$ are Lipschitz continuous with Lipschitz constant $\lambda_i$, $i\in\mathcal{S}$.
\end{hyp}
\begin{hyp}\label{hyp:infinite}
The state space $\mathcal{S}$ is countable, the quantities $\|A\|_{\mathcal{S}}$, $\|C\|_{\mathcal{S}}$, $\|D\|_{\mathcal{S}}$,  $\|DD'\|_{\mathcal{S}}$ and $\|E\|_{\mathcal{S}}$ are finite. The cumulative distribution functions of the sojourn times $F_i$ are Lipschitz continuous with Lipschitz constant $\lambda_i$, $i\in\mathcal{S}$ and
\begin{equation*}
\overline{\lambda}=\sup_{i\in\mathcal{S}}\{\lambda_{i}\}<\infty.
\end{equation*}
\end{hyp}
Note that the extra assumptions in the infinite case hold true automatically in the finite case, 
and that the Lipschitz assumptions hold true automatically for MJLS (i.e., when 
the distributions of $F_i$ are exponential).

We address the filtering problem of estimating the value of $x(t)$ given the observations $\{y(s),\theta(s), 0\le s\le t\}$ for $0\leq t\leq T$. 
It is well-known that the KBF is the optimal estimator because 
the problem is equivalent to estimating the state of a linear time-varying system (with no jumps), 
taking into account that the past values of $\theta$ are available.
The KBF satisfies the following equation 
\begin{equation*}
d\hat{x}_{KB}(t)\!=\!A_{\theta(t)}\hat{x}_{KB}(t)dt+K_{KB}(t)(dy(t)-C_{\theta(t)}\hat{x}_{KB}(t)dt),
\end{equation*}
for $0\leq t\leq T$, with initial condition $\hat{x}_{KB}(0)=\mathbb{E}[x(0)]$ and gain matrix 
\begin{equation}\label{eq Ricc gain}
K_{KB}(t)=P_{KB}(t)C_{\theta(t)}'(D_{\theta(t)}D_{\theta(t)}')^{-1},
\end{equation}
for $0\leq t\leq T$, where $P_{KB}(t)$ is an $n_1\times n_1$ matrix-valued process satisfying the Riccati matrix differential equation
\begin{equation}
\label{eq Ric}
\left\{\begin{array}{rcl}
dP_{KB}(t)&=&R(P_{KB}(t),\theta(t))dt,\\ 
P_{KB}(0)&=&Var(x(0)),
\end{array}\right.
\end{equation}
for $0\leq t\leq T$, where $R:\dR^{n_1\times n_1}\times \mathcal S \rightarrow \dR^{n_1\times n_1}$
is defined for any $M\in\dR^{n_1\times n_1}$ and $i\in\mathcal{S}$ by
\begin{equation}\label{eq:def:R}
R(M,i)=A_iM+MA_i'+E_iE_i'-MC_i'(D_iD_i')^{-1}C_iM.
\end{equation}

 It is usually not possible to pre-compute a solution for this system 
(prior to the observation of $\theta(s)$, $0\le s\le t$).
Moreover, to solve it in real time after observing $\theta$ would require 
instantaneous computation of $P(t)$; one can obtain a delayed solution 
 $P(t-\delta)$ where $\delta$ is the time required to solve the system, however 
using this solution as if it was the actual $P(t)$ in the filter
may bring considerable error to the obtained estimate depending on $\delta$ and on the 
system parameters (e.g., many jumps may occur between $t-\delta$ and $t$).

The aim of this paper is to propose a new filter based on suitably chosen pre-computed solutions of Eq.~(\ref{eq Ric}) 
under the finiteness assumption~\ref{hyp:finite}
and to show convergence of our estimate to the optimal KBF when the number of discretization points goes to infinity
under the more general countable assumption~\ref{hyp:infinite}.
We also compare its performance with the Fragoso-Costa LMMSE filter \cite{FC10} on a real-world application.
%
\section{Approximate Kalman--Bucy filter}
\label{sec-KBF}
The estimator is constructed as follows. 
We first select an optimized finite set of typical possible trajectories of $\{\theta(t)$, $0\le t\le T\}$ by 
discretizing the {semi-}Markov chain and for each such trajectory we solve Eqs.~(\ref{eq Ric}), \eqref{eq Ricc gain} 
and store the results. 
In real time, the estimate is obtained by looking up the pre-computed solutions 
and selecting the suitable gain given the current value of $\theta(t)$.
%
\subsection{Discretization of the {semi-}Markov chain}
The approach relies on the construction of optimized typical trajectories of the {semi-}Markov chain $\{\theta(t), 0\leq t\leq T\}$. 
First we need to rewrite this {semi-}Markov chain in terms of its jump times and post-jump locations.
Let $T_0=0$ and $T_k$ be the $k$-th jump time of $\{\theta(t), 0\leq t\leq T\}$ for $k\geq 1$,
\begin{equation*}
T_{k}=\inf\{t\geq T_{k-1};\ \theta(t)\neq\theta(T_{k-1})\}.
\end{equation*}
For $k\geq 0$ let $Z_k=\theta(T_k)$ be the post-jump locations of the chain. 
Let $S_0=0$ and for $k\geq 1$, $S_k=T_k-T_{k-1}$ be the inter-arrival times of 
the Markov process $\{\theta(t), 0\leq t\leq T\}$. Using this notation, $\theta(t)$ can be rewritten as
\begin{equation}\label{eq:theta}
\theta(t)=\sum_{k=0}^\infty Z_k\ind{T_k\leq t<T_{k+1}}=\sum_{k=0}^\infty Z_k\ind{0\leq t-T_k<S_{k+1}}.
\end{equation}
Under the finiteness assumption~\ref{hyp:finite},
as the state space $\mathcal{S}$ of $\{\theta(t), 0\leq t\leq T\}$ (and hence of $\{Z_k\}$) is finite, to obtain a fully discretized approximation of $\{\theta(t), 0\leq t\leq T\}$ one only needs to discretize the inter-arrival times $\{S_k\}$ on a finite state space. One thus constructs a finite set of typical possible trajectories of $\{\theta(t), 0\leq t\leq T\}$ up to a given jump time horizon $T_n$ selected such that $T_n\geq T$ with high enough probability. 

To discretize the inter-arrival times $\{S_k\}$, we choose a quantization approach that has been recently developed in numerical probability.
Its main advantage is that the discretization is optimal in some way explained below.
There exists an extensive literature on quantization methods for random variables and processes. The interested reader may for instance, consult the following works \cite{bally03,gray98,pages98} and references therein.
Consider $X$ an $\mathbb{R}^m$-valued random variable such that $\dE[| X |^2] < \infty$ 
and $\nu$ a fixed integer; 
the optimal $L^{2}$-quantization of the random variable $X$ consists in finding the best possible $L^{2}$-approximation of $X$ by a random vector $\widehat{X}$ taking at most $\nu$ different values, which can be carried out in two steps.
First, find a finite weighted grid $\Gamma\subset \mathbb{R}^{m}$ with $\Gamma= \{\gamma^{1},\ldots,\gamma^{\nu}\}$. 
Second, set $\widehat{X}=\widehat{X}^{\Gamma}$ where $\widehat{X}^{\Gamma}=proj_{\Gamma}(X)$ with $proj_{\Gamma}$ denoting the closest neighbor projection on $\Gamma$.
The asymptotic properties of the $L^{2}$-quantization are given in e.g. \cite{pages98}.
\begin{theorem}\label{th:quantize}
If $\mathbb{E}[|X|^{2+\epsilon}]<+\infty$ for some $\epsilon>0$ then one has
\begin{eqnarray*}
\lim_{\nu\rightarrow \infty} \nu^{1/m} \min_{|\Gamma|\leq \nu} \dE[| X-\widehat{X}^{\Gamma}|^{2}]^{1/2}& =& C,
\end{eqnarray*}
for some contant $C$ depending only on $m$ and the law of $X$ and where $|\Gamma|$ denote the cardinality of $\Gamma$.
\end{theorem}

Therefore the $L^2$ norm of the difference between $X$ and its quantized approximation $\wh{X}$ goes to zero with rate $\nu^{-1/m}$ as the number of points $\nu$ in the quantization grid goes to infinity. The competitive learning vector quantization algorithm (CLVQ) provides the optimal grid based on a random simulator of the law of $X$ and a stochastic gradient method.

In the following, we will denote by $\wh{S}_k$ the quantized approximation of the random variable $S_k$ and $\wh{T}_k=\wh{S}_1+\cdots+\wh{S}_k$ for all $k$. 
%
\subsection{Pre-computation of a family of solutions to Riccati equation}
We start by rewriting the Riccati equation~(\ref{eq Ric}) in order to have a similar 
expression to Eq.~(\ref{eq:theta}). 
As operator $R$ does not depend on time, the solution $\{P(t), 0\leq t\leq T\}$ to Eq.~(\ref{eq Ric}) corresponding to a given trajectory $\{\theta(t), 0\leq t\leq T\}$ can be rewritten as
\begin{equation*}
P(t)=\sum_{k=0}^\infty P_k(t-T_k)\ind{0\leq t-T_k<S_{k+1}},
\end{equation*}
for $0\leq t\leq T$, where $\{P_0(t), 0\leq t\leq T\}$ is the solution of the system
\begin{equation*}
\left\{\begin{array}{rcl}
d{P}_0(t)&=&R(P_0(t),Z_0)dt,\\
P_0(0)&=&p_0,
\end{array}\right.
\end{equation*}
for $0\leq t\leq T$, with $p_0=Var(x(0))$, 
and for $k\geq 1$, $\{P_k(t), 0\leq t\leq T\}$ is recursively defined as the solution of 
\begin{equation*}
\left\{\begin{array}{rcl}
d{P}_k(t)&=&R(P_k(t),Z_k)dt,\\
P_k(0)&=&P_{k-1}(S_k).
\end{array}\right.
\end{equation*}
Given the quantized approximation $\{\wh{S}_k\}$ of the sequence $\{{S}_k\}$, we propose the following approximations $\{\wh{P}_k(t), 0\leq t\leq T\}$ of $\{P_k(t), 0\leq t\leq T\}$ for all $k$.
First, $\{\wh{P}_0(t), 0\leq t\leq T\}$ is the solution of
\begin{equation*}
\left\{\begin{array}{rcl}
d{\wh{P}}_0(t)&=&R(\wh{P}_0(t),Z_0)dt,\\
\wh{P}_0(0)&=&p_0,
\end{array}\right.
\end{equation*}
and for $k\geq 1$, $\{\wh{P}_k(t), 0\leq t\leq T\}$ is recursively defined as the solution of 
\begin{equation*}
\left\{\begin{array}{rcl}
d{\wh{P}}_k(t)&=&R(\wh{P}_k(t),Z_k)dt,\\
\wh{P}_k(0)&=&\wh{P}_{k-1}(\wh{S}_k).
\end{array}\right.
\end{equation*}
Hence $P_k$ and $\wh{P}_k$ are defined with the same dynamics, the same horizon $T$, but different starting values, and all the $\wh{P}_k$ can be computed off-line for each of the finitely many possible values of $(Z_k,\wh{S}_k)$ 
(under the finiteness assumption~\ref{hyp:finite} and for a finite number of jumps) 
and stored.
%
\subsection{On line approximation}
We suppose that on-line computations are made on a regular time grid with constant step $\delta t$. Note that in most applications $\delta t$ is small compared to the time $\delta$ of instantaneous computation of $P(t)$. The state of the {semi-}Markov chain $\{\theta(t), 0\leq t\leq T\}$ is observed, 
but the jumps can only be considered, in the filter operation,
at the next point in the time grid. Set $\wt{T}_0=0$, and for $k\geq 1$ define $\wt{T}_k$ as
\begin{equation*}
\wt{T}_{k}=\inf\{j;\  T_k< j\delta t\}\delta t,
\end{equation*}
hence $\wt{T}_{k}$ is the effective time at which the $k$-th jump is taken into account. One has $\wt{T}_{k}>T_k$ and the difference between $\wt{T}_{k}$ and $T_k$ is at most $\delta t$. We also set $\wt{S}_k=\wt{T}_k-\wt{T}_{k-1}$ for $k\geq 1$.
Now we construct our approximation $\{\wt P(t), 0\leq t\leq T\}$ of $\{P(t), 0\leq t\leq T\}$ as follows
\begin{equation*}
\wt{P}(t)=\sum_{k=0}^\infty \wh{P}_k(t-\wt{T}_k)\ind{0\leq t-\wt{T}_k<\wt{S}_{k+1}}\ind{t\leq T}.
\end{equation*}
Thus we just select the appropriate pre-computed solutions and paste them at the approximate jumps times $\{\wt{T}_{k}\}$, which can be done on-line. The approximate gain matrices are simply defined by
\begin{equation*}
\wt{K}(t)=\wt{P}(t)C_{\theta(t)}'(D_{\theta(t)}D_{\theta(t)}')^{-1},
\end{equation*}
and the estimated trajectory satisfies
\begin{equation*}
d\wt{x}(t)=A_{\theta(t)}\wt{x}(t)dt+\wt{K}(t)(dy(t)-C_{\theta(t)}\wt{x}(t)dt),
\end{equation*}
for $0\leq t\leq T$, with initial condition $\wt{x}(0)=\mathbb{E}[x(0)]$.
%
\section{Convergence of the approximation procedure}
\label{sec-convergence}
The investigation of the convergence of our approximation scheme
under the general assumption~\ref{hyp:infinite},
is made in several
steps again. The first one is the evaluation of the error between $P(t)$ and $\wt{P}(t)$ up to the time horizon $T$ and requires some Lipschitz regularity assumptions on the solution of Riccati equations. First, we establish these regularity properties. Then we derive the error between $P$ and $\wt{P}$, and finally we evaluate the error between the real KBF filter $\wh{x}_{KB}$ and its quantized approximation $\wt{x}$.
%
\subsection{Regularity of the solutions of Riccati equations}
For all $t\geq 0$, suitable matrix $p\in\mathcal{C}(n_1)$ and $i\in\mathcal{S}$ denote by $\phi_i(p,t)$ the solution at time $t$ of the following Riccati equation starting from $p$ at time $0$,
\begin{equation*}\label{eqi Riccati}
\left\{\begin{array}{rcl}
d{P}(t)&=&R(P(t),i)dt,\\
P(0)&=&p,
\end{array}\right.
\end{equation*}
for $t\geq 0$. 
We start with a boundedness result.
\begin{lemma}\label{lem00:Lip}
Under Assumption~\ref{hyp:infinite}, 
for all $\bar p_0\in\mathcal{C}(n_1)$, there exist a matrix $\bar p_1\in\mathcal{C}(n_1)$ such that $\bar p_1\geq \bar p_0$ and for $p\leq \bar p_0$, $i\in\mathcal{S}$ and times $0\leq t\leq T$, one has
$\phi_i(p,t)\leq \bar{p}_1$.
\end{lemma}
\textit{Proof.}  
The Riccati equation can be rearranged in the following form 
\begin{eqnarray*}
\frac{d P(t)}{dt} &=&   A_{aux}(t)P(t)  + P(t) A_{aux}(t)' + E_iE_i' \\
&&+ K_i(t)D_iD_i'K_i(t)',
\end{eqnarray*}
where $K_i(t)= P(t)C_i'(D_iD_i')^{-1}$ and $A_{aux}(t)=A_i - K_i(t)C_i$. 
For any matrix $L$ with suitable dimensions, from the optimality of the KBF we have that  
$\phi_i(p,t)\leq \phi_L(p,t)$ where $\phi_L(p,t)$ is the covariance 
of a linear state observer with gain $L$, so that  $\phi_L(p,t)$ is the solution of 
\begin{eqnarray*}
\frac{d P(t)}{dt} &= &  (A_i-LC_i)(t)P(t)  + P(t) (A_i-LC_i)' \\
&&+ E_iE_i' + LD_iD_i'L',\\
P(0)&=&p. 
\end{eqnarray*}
In particular, we can set $L=0$, and $\phi_L(p,t)$ is now the solution of the linear differential equation
\begin{equation}\label{eq:covariance:of:trivial:filter}
\frac{d P(t)}{dt} =   A_iP(t)  + P(t) A_i' + E_iE_i',  \qquad P(0)=p,
\end{equation}
which can be expressed in the form $\phi_L(p,t)=\Phi_1(t)p+\Phi_2(t)$ where 
$\Phi_1\leq \beta e^{\alpha {\|A_i\|}t} \|p\| I$ and 
$\Phi_2\leq \int_0^t \beta e^{\alpha {\|A_i\|}\tau} \|E_iE_i\| I d\tau $ 
for some scalars $\alpha,\beta$ that do not depend on $p,i$.   
Set $\bar p_1=\beta e^{\alpha T{\|A\|_{\mathcal{S}}}}(\|\bar p_0\|p_0+T \|E\|_{\mathcal{S}}^2 I )$, thus completing the proof. 
\qquad\hspace{\stretch{1}}$ \Box$

\begin{theorem}\label{lem-1-Ric}
Under Assumption~\ref{hyp:infinite}, for each $\wt p\in\mathcal{C}(n_1)$ 
there exist $\ell,\eta>0$ such that for all $i\in\mathcal{S}$  
and $0\leq t,\wh{t}\leq T$ and $p,\wh{p}\leq \wt p$ one has
\begin{equation*}
\|\phi_i(p,t)-\phi_i(\wh p,\wh t)\|\leq \ell |t-\wh{t}|+\eta\|p-\wh{p}\|.
\end{equation*}
\end{theorem}
\textit{Proof.}  
It follows directly from the definition of $R$ in Eq.~\eqref{eq:def:R} that 
one has
\begin{eqnarray*}
\lefteqn{\frac{d\phi_i(p,t)-d\phi_i(\wh{p},t)}{dt}}\\
&=& A_i \phi_i(p,t) + \phi_i(p,t) A_i'+ E_iE_i' \\
&&- \phi_i(p,t) C_i'(D_iD_i')^{-1}C_i \phi_i(p,t)\\
&&- \big( A_i \phi_i(\wh{p},t) + \phi_i(\wh{p},t) A_i'+ E_iE_i' \\
&&- \phi_i(\wh{p},t) C_i'(D_iD_i')^{-1}C_i \phi_i(\wh{p},t) \big)\\ 
&=& A_i (\phi_i(p,t)-\phi_i(\wh{p},t)) + (\phi_i(p,t)-\phi_i(\wh{p},t)) A_i'\\
&&- \phi_i(\wh{p},t) C_i'(D_iD_i')^{-1}C_i (\phi_i(p,t)-\phi_i(\wh{p},t))\\
&&- (\phi_i(p,t)-\phi_i(\wh{p},t)) C_i'(D_iD_i')^{-1}C_i \phi_i(\wh{p},t)\\
&& - (\phi_i(p,t)-\phi_i(\wh{p},t)) C_i'(D_iD_i')^{-1}C_i\\
&&\times (\phi_i(p,t)-\phi_i(\wh{p},t))  \\  
&=& (A_i - \phi_i(\wh{p},t) C_i'(D_iD_i')^{-1}C_i) (\phi_i(p,t)-\phi_i(\wh{p},t)) \\
&&+ (\phi_i(p,t)-\phi_i(\wh{p},t)) (A_i'-C_i'(D_iD_i')^{-1}C_i \phi_i(\wh{p},t))\\
&& - (\phi_i(p,t)-\phi_i(\wh{p},t)) C_i'(D_iD_i')^{-1}C_i \\
&&\times(\phi_i(p,t)-\phi_i(\wh{p},t)),
\end{eqnarray*}
or, by denoting $X(t)=\phi_i(p,t)-\phi_i(\wh{p},t)$, one has $X(0)=p-\wh{p}$ and 
\begin{eqnarray}
\frac{d X(t)}{dt} &=&   A_{aux}(t)X(t)  + X(t) A_{aux}(t)' \nonumber\\
&&- X(t) C_i'(D_iD_i')^{-1}C_i X(t), \label{eq:aux:X}
\end{eqnarray}
where we write $A_{aux}(t)=(A_i - \phi_i(\wh{p},t) C_i'(D_iD_i')^{-1}C_i)$ 
for ease of notation. 
By setting $Y(0)=\|p-\wh{p}\|I\geq X(0)$ and using the order preserving property of 
the Riccati equation \eqref{eq:aux:X} it follows that $\{Y(t), 0\leq t\leq T\}$ defined as the solution of 
\begin{eqnarray}
\frac{d Y(t)}{dt} &=&  A_{aux}(t) Y(t) + Y(t) A_{aux}(t)' \nonumber\\
&& - Y(t) C_i'(D_iD_i')^{-1}C_i Y(t), \label{eq:aux:Y}
\end{eqnarray}
satisfies $Y(t)\geq X(t)$ for all $t\geq 0$.
The process $\{Y(t), 0\leq t\leq T\}$ can be interpreted as the error covariance of a filtering problem%
\footnote{Note that this does not hold true for the process $\{X(t), 0\leq t\leq T\}$ as it may not be positive semidefinite.}, 
more precisely the covariance of the error $\wh{x}_{aux}(t)-x_{aux}(t)$ where $\{\wh{x}_{aux}(t), 0\leq t\leq T\}$ satisfies
\begin{equation*}
d \wh{x}_{aux} = A_{aux}(t) \wh{x}_{aux} dt + K(t) (dy - C_{aux} \wh{x}_{aux} dt),
\end{equation*}
with $A_{aux}$ defined above, 
$C_{aux}=(C_i'(D_iD_i')^{-1}C_i)^{1/2}$, $\{K(t), 0\leq t\leq T\}$ is the Kalman gain, and
\begin{equation*}
\left\{\begin{array}{rcl}
d x_{aux}(t) &=& A_{aux}(t) x_{aux}(t)dt ,\\
d y_{aux}(t) &=& C_{aux} x_{aux}(t)dt + dv_{aux}(t),
\end{array}\right.
\end{equation*}
where $\{v_{aux}(t), 0\leq t\leq T\}$ is a standard Wiener process with incremental covariance $Idt$, 
and $x_{aux}(0)$ is a Gaussian random variable with covariance $p-\wh{p}$.
Now, if we replace $K$ with the (suboptimal) gain $L=0$
we obtain a larger error covariance $Y_{L} (t)\geq Y(t)$. 
With the trivial gain $L=0$ we also have 
\begin{equation*}
d\wh{x}_{aux}- dx_{aux} = A_{aux}(t)(\wh{x}_{aux}-x_{aux}) dt,
\end{equation*}
so that direct calculation yields 
\begin{equation}\label{eq:aux:YL}
\frac{d Y_L(t)}{dt} = A_{aux}(t) Y_L(t) + Y_L(t) A_{aux}(t)',
\end{equation}
with $Y_L(0)=\|p-\wh{p}\|I$. Recall that $\wh p\leq \wt p$ by hypothesis, 
so that from Lemma \ref{lem00:Lip} we get an uniform bound $\bar p_1$ for 
$\phi_i(\wh p,t)$, which in turn yields
that  $\|A_{aux}\|_{\mathcal{S}}$ is bounded in the time interval $0\leq t\leq T$ 
and for all $\wh p\leq \wt p$. This allows to write
\begin{equation*}
Y(t) \leq \ell_1 \|p-\wh{p}\| I,\qquad 0\leq t\leq T,
\end{equation*}
for some $\ell_1\geq 0$ (uniform on $t$, $p$, $\wh p$ and $i$). 
Gathering some of the above inequalities together, one gets
\begin{equation}\label{eq-aux-main-eval01}
\phi_i(p,t)-\phi_i(\wh{p},t)=X(t)\leq Y(t)\leq Y_L(t)\leq \ell_1 \|p-\wh{p}\| I,
\end{equation}
$0\leq t\leq T$. Similarly as above, one can obtain 
\begin{equation}\label{eq-aux-main-eval02}
\phi_i(\wh p,t)-\phi_i(p,t)\leq \ell_2 \|p-\wh{p}\| I,\qquad 0\leq t\leq T,
\end{equation}
where, again, $\ell_2$ is uniform on $t$, $p$, $\wh p$ and $i$.
Eqs.~\eqref{eq-aux-main-eval01}, \eqref{eq-aux-main-eval02} 
and the fact that 
$\phi_i(\wh p,t)-\phi_i(p,t)$ is symmetric lead to
\begin{eqnarray*}
-\max(\ell_1,\ell_2)&\leq& \lambda_{\min}(\phi_i(\wh p,t)-\phi_i(p,t)),\\
\lambda_{\min}(\phi_i(\wh p,t)-\phi_i(p,t))&\leq &\lambda_{\max}(\phi_i(\wh p,t)-\phi_i(p,t)) \\
&\leq &\max(\ell_1,\ell_2).
\end{eqnarray*}
Hence, one has 
\begin{equation*}
\|\phi_i(\wh p,t)-\phi_i(p,t)\|\leq \max(\ell_1,\ell_2) \|p-\wh{p}\|,
\end{equation*}
completing the first part of the proof. 

For the second part, similarly to the proof of the preceding lemma, 
we have that $\phi_i(p,t)$ is bounded from above by $X(t)$ 
the solution of the linear differential equation 
Eq.\eqref{eq:covariance:of:trivial:filter} with initial condition $X(0)=p$,
and it is then simple to find scalars $\eta_1,\eta_2>0$ irrespective of 
$i$ such that, for the entire time interval $0\leq t\leq T$,
\begin{equation*}
\|X(t)-p\|_{\mathcal{S}} \leq \|\Phi_1(t)\|_{\mathcal{S}} + \|(\Phi_2(t)-I)p\|_{\mathcal{S}}  
\leq \eta_1 t +  \eta_2  t \|p\|.
\end{equation*}
Hence, one has
\begin{equation}\label{eq:maj phi}
\phi_i(p,t)-p\leq X(t)-p \leq \|X(t)-p\|_{\mathcal{S}} I \leq (\eta_1 t +  \eta_2  t \|p\|) I,
\end{equation}
for all $t\geq 0$, leading to
\begin{equation*}
\|\phi_i(p,t)-p\| \leq \eta_1 t +  \eta_2  t \|p\|.
\end{equation*}
As $p\leq \wt p$ by hypothesis, we have $\|p\|\leq n_1\|\wt p\|$ and 
it follows immediately from the above inequality that
\begin{equation}\label{eq:aux:XLips2}
\|\phi_i(p,t)-p\| \leq (\eta_1 +  \eta_2 n_1 \|\wt p\|)  t .
\end{equation}
As operator $R$ does not depend on time, 
we have $\phi(p,t_1+t_2)=\phi(\phi(p,t_1),t_2)$, $t_1,t_2\geq 0$,
and defining $\bar p=\phi(p,t_1)$, one has
\begin{equation*}
\|\phi_i(p,t_1+t_2)-\phi_i(p,t_1)\| = \|\phi_i(\bar p,t_2)-\bar p\| 
\end{equation*}
and Eq.~\eqref{eq:aux:XLips2} allows to write
\begin{equation*}
\|\phi_i(p,t_1+t_2)-\phi_i(p,t_1)\|  \leq (\eta_1 +  \eta_2 n_1 \|\wt p\|) t_2.
\end{equation*}
The result then follows by setting $t_1=\wh t$ and $t_2=t-\wh t$ if $t> \wh t$
or with $t_1=t$ and $t_2=\wh t - t$ otherwise. 
\qquad\hspace{\stretch{1}}$ \Box$
%
\subsection{Error derivation for gain matrices}
We proceed in three steps. The first one is to study the error between $P_k(t)$ and $\wh{P}_k(t)$, the second step is to study the error between $P(t)$ and $\wt{P}(t)$ and the last step is to compare the gain matrices $K_{KB}(t)$ and $\wt{K}(t)$, for $0\leq t\leq T$. We start with a preliminary important result that will enable us to use Theorem~\ref{lem-1-Ric}
in all the sequel.
\begin{lemma}\label{lem0:Lip}
Under Assumption~\ref{hyp:infinite}, there exist a matrix $\bar p\in\mathcal{C}(n_1)$ such that for all integers $0\leq k\leq n$ and times $0\leq t\leq T$, one has 
\begin{equation*}
P_k(t)\leq \bar p,\qquad \wh{P}_k(t)\leq \bar p.
\end{equation*}
\end{lemma}
\textit{Proof.}  We prove the result by induction on $k$. For $k=0$, one has $p_0\in\mathcal{C}(n_1)$ and $P_0(t)=\wh{P}_0(t)=\phi_{Z_0}(p_0,t)$ for all $t\leq T$. Lemma~\ref{lem00:Lip} thus yields the existence of a matrix $\bar p_0\in\mathcal{C}(n_1)$ such that $P_0(t)\leq \bar p_0$ for all $t\leq T$.
Suppose that for a given $k\leq n-1$, there exists a matrix $\bar p_k\in\mathcal{C}(n_1)$ such that $P_k(t)\leq \bar p_k$ and $\wh{P}_k(t)\leq \bar p_k$ for all $t\leq T$. Then in particular, if $S_k\leq T$ and $\wh{S}_k\leq T$, one has $P_{k+1}(0)=P_k(S_k)\leq \bar p_k$ and $\wh{P}_{k+1}(0)=\wh{P}_k(\wh{S}_k)\leq \bar p_k$. Hence, Lemma~\ref{lem00:Lip} gives the existence of a matrix $\bar p_{k+1}\in\mathcal{C}(n_1)$ such that $P_{k+1}(t)\leq \bar p_{k+1}$ and $\wh{P}_{k+1}(t)\leq \bar p_{k+1}$ for all $t\leq T$. One thus obtains an increasing sequence $(p_k)$ of matrices in $\mathcal{C}(n_1)$ and the result is obtained by setting $\bar p=\bar p_n$.
\qquad \hspace{\stretch{1}}$ \Box$

In the following, for  $\bar p$ given by Lemma~\ref{lem0:Lip} we set $\tilde p= \bar p$
in Theorem~\ref{lem-1-Ric}
and denote by $\bar \ell$ and $\bar \eta$ the 
corresponding Lipschitz constants.
We now turn to the investigation of the error between the processes $P_k(t)$ and $\wh{P}_k(t)$.
\begin{lemma}\label{lem1:Lip}
Under Assumption~\ref{hyp:infinite}, for all integers $0\leq k\leq n$ and times $0\leq t\leq T$, one has
\begin{equation*}
\|P_k(t)-\wh{P}_k(t)\|\leq\bar\ell\|P_{k-1}(S_k)- \wh{P}_{k-1}(\wh{S}_k)\|.
\end{equation*}
\end{lemma}
\textit{Proof.}  One has $P_k(t)=\phi_{Z_k}(P_{k-1}(S_k),t)$ and $\wh{P}_k(t)=\phi_{Z_k}(\wh{P}_{k-1}(\wh{S}_k),t)$. Hence, Lemma~\ref{lem0:Lip} and Theorem~\ref{lem-1-Ric}
yield
\begin{eqnarray*}
\lefteqn{\|P_k(t)-\wh{P}_k(t)\|}\\
&=&\|\phi_{Z_k}(P_{k-1}(S_k),t)-\phi_{Z_k}(\wh{P}_{k-1}(\wh{S}_k),t)\|\\
&\leq&\bar\ell\|P_{k-1}(S_k)- \wh{P}_{k-1}(\wh{S}_k)\|,
\end{eqnarray*}
if $S_k,\wh{S}_k\leq T$, hence the result.
\qquad \hspace{\stretch{1}}$ \Box$
\begin{lemma}\label{lem2:Lip}
Under Assumption~\ref{hyp:infinite}, for all integers $0\leq k\leq n$ satisfying $S_k,\wh{S}_k\leq T$, one has
\begin{equation*}
\|P_{k}(S_{k+1})- \wh{P}_{k}(\wh{S}_{k+1})\|\leq \sum_{j=0}^{k}\bar\ell^{k-j}\bar\eta|S_{j+1}-\wh{S}_{j+1}|.
\end{equation*}
\end{lemma}
\textit{Proof.}  By definition, one has $P_k(S_{k+1})=\phi_{Z_k}(P_{k-1}(S_k),S_{k+1})$ and $\wh{P}_k(t)=\phi_{Z_k}(\wh{P}_{k-1}(\wh{S}_k),\wh{S}_{k+1})$. Hence as above, one has
\begin{eqnarray*}
\lefteqn{\|P_{k}(S_{k+1})- \wh{P}_{k}(\wh{S}_{k+1})\|}\\
&=&\|\phi_{Z_k}(P_{k-1}(S_k),S_{k+1})-\phi_{Z_k}(\wh{P}_{k-1}(\wh{S}_k),\wh{S}_{k+1})\|\\
&\leq&\bar\ell\|P_{k-1}(S_k)- \wh{P}_{k-1}(\wh{S}_k)\|+\bar\eta|S_{k+1}-\wh{S}_{k+1}|.
\end{eqnarray*}
Then notice that one also has
\begin{eqnarray*}
\lefteqn{\|P_0(S_1)-\wh{P}_0(\wh{S}_1)\|}\\
&=&\|\phi_{Z_0}(p_0,S_1)-\phi_{Z_0}(p_0,\wh{S}_1)\|\leq\bar\eta|S_1-\wh{S}_1|,
\end{eqnarray*}
and the result  is obtained by recursion.
\qquad \hspace{\stretch{1}}$ \Box$

We can now turn to the error between the processes $P(t)$ and $\wt{P}(t)$.
\begin{theorem}\label{th:Lip}
Under Assumption~\ref{hyp:infinite}, for all $0\leq t< T\wedge T_{n+1}$, one has
\begin{eqnarray*}
\lefteqn{\dE[\|P(t)-\wt{P}(t)\|^2\ind{0\leq t\leq T\wedge T_{n+1}}]^{1/2}}\\
&\leq&\sum_{j=0}^{n-1}\bar\ell^{n-j}\bar\eta\dE[|S_{j+1}-\wh{S}_{j+1}|^2]^{1/2}\\
&&+\bar\eta \delta t+n\|\bar p\|(\overline{\lambda}\delta t)^{1/2},
\end{eqnarray*}
where $\bar p$ is defined in Lemma~\ref{lem0:Lip}.
\end{theorem}
\begin{rem}\label{rem-byproduct}
Note that the above result is very general. Indeed, we do not use in its proof that $\wh{S}_k$ is the quantized approximation of $S_k$. We have established that, given a {semi-}Markov chain $\{\theta(t), 0\leq t\leq T\}$ and a process $\{\wh{\theta}(t),0\leq t\leq T\}$ obtained by a perturbation of the jump times of $\{\theta(t), 0\leq t\leq T\}$, the two solutions of the Riccati equations driven by these two processes respectively are not far away from each other, as long as the real and perturbed jump times are not far away from each other. 
We allow two kinds of perturbations, a random one, given by the replacement of $S_k$ by $\wh{S}_k$ 
and a deterministic one given by $\delta t$ corresponding to a delay in the jumps. 
In the case of non-instantaneous observation of $\theta(t)$ (i.e., imperfect observation $\wt{S}_k$ of 
$S_k$), the difference $\dE[|\wt{S}_{j+1}-\wh{S}_{j+1}|^2]$ may not converge to zero but is still a valid upper bound for the approximation error of the Riccati solution and can reasonably be supposed small enough.
Note also that the result is still valid for any $L^q$ norm instead of the $L^2$ norm as the initial value of the Riccati solution is deterministic, as long as the distributions $F_i$ have moments of order greater than $q$.
\end{rem}
\textit{Proof.}  By definition, one has for all $0\leq t< T\wedge T_{n+1}$
\begin{eqnarray*}
\lefteqn{P(t)-\wt{P}(t)}\\
&=&\sum_{k=0}^nP_k(t-T_k)\ind{0\leq t-T_k<S_{k+1}}\\
&&-\wh{P}_k(t-\wt{T}_k)\ind{0\leq t-\wt{T}_k<\wt{S}_{k+1}}\\
&=&\sum_{k=0}^n\big(P_k(t-T_k)-\wh{P}_k(t-{T}_k)\big)\ind{0\leq t-T_k<S_{k+1}}\\
&&+\sum_{k=0}^n\big(\wh{P}_k(t-{T}_k-\wh{P}_k(t-\wt{T}_k)\big)\ind{0\leq t-T_k<S_{k+1}}\\\
&&+\sum_{k=0}^n\!\wh{P}_k(t-\wt{T}_k)(\ind{0\leq t-T_k<S_{k+1}}\!-\!\ind{0\leq t-\wt{T}_k<\wt{S}_{k+1}})\\
&=&\epsilon_1(t)+\epsilon_2(t)+\epsilon_3(t).
\end{eqnarray*}
From Lemmas \ref{lem1:Lip} and \ref{lem2:Lip}, the first term $\epsilon_1$ can be bounded by
\begin{eqnarray*}
\lefteqn{\|\epsilon_1(t)\|}\\
&\leq&\big\|\sum_{k=0}^n \big(P_k(t-T_k)-\wh{P}_k(t-T_k)\big)\ind{0\leq t-T_k<S_{k+1}}\big\|\\
&\leq& \sum_{k=0}^n \|P_k(t-T_k)-\wh{P}_k(t-T_k)\|\ind{0\leq t-T_k<S_{k+1}}\\
&\leq& \sum_{k=0}^n \ell\|P_{k-1}(S_k)- \wh{P}_{k-1}(\wh{S}_k)\|\ind{0\leq t-T_k<S_{k+1}}\\
&\leq& \sum_{k=0}^n \sum_{j=0}^{k-1}\bar\ell^{k-j}\bar\eta|S_{j+1}-\wh{S}_{j+1}|\ind{T_k\leq t<T_{k+1}}\\
&\leq& \sum_{j=0}^{n-1}\bar\ell^{n-j}\bar\eta|S_{j+1}-\wh{S}_{j+1}|.
\end{eqnarray*}
The second term $\epsilon_2$ is bounded by Lemma~\ref{lem0:Lip} and Theorem~\ref{lem-1-Ric}
as follows
\begin{eqnarray*}
\lefteqn{\|\epsilon_2(t)\|}\\
&\leq&\big\|\sum_{k=0}^n \big(\wh{P}_k(t-{T}_k-\wh{P}_k(t-\wt{T}_k)\big)\ind{0\leq t-T_k<S_{k+1}}\big\|\\
&\leq& \sum_{k=0}^n \|\wh{P}_k(t-T_k)-\wh{P}_k(t-\wt{T}_k)\|\ind{0\leq t-T_k<S_{k+1}}\\
&\leq& \sum_{k=0}^n\bar\eta |T_k-\wt{T}_k|\ind{0\leq t-T_k<S_{k+1}}\\
&\leq& \bar\eta\delta t,
\end{eqnarray*}
using the fact that the difference between $T_k$ and $\wt{T}_k$ is less than $\delta t$ by construction. Finally, the last term $\epsilon_3$ is bounded by using Lemma~\ref{lem0:Lip} and the fact that $0\leq {T}_k\leq \wt{T}_k$ for all $k$. Indeed, one has
\begin{eqnarray*}
\lefteqn{\dE[\|\epsilon_3(t)\|^2]^{1/2}}\\
&\leq&\dE\big[\big\|\sum_{k=0}^n\wh{P}_k(t-\wt{T}_k)(\ind{0\leq t-T_k<S_{k+1}}\\
&&-\ind{0\leq t-\wt{T}_k<\wt{S}_{k+1}})\big\|^2\big]^{1/2}\\
&\leq&\|\bar p\| \sum_{k=0}^n \dE[|\ind{0\leq t-T_k<S_{k+1}}-\ind{0\leq t-\wt{T}_k<\wt{S}_{k+1}}|^2]^{1/2}\\
&\leq&\|\bar p\| \sum_{k=0}^n \dP(t-\delta t\leq T_k\leq t)^{1/2}\\
&\leq& n\|\bar p\|\sum_{i\in\mathcal{S}}\big({\lambda_i \delta t}\big)^{1/2}\dP(Z_k=i)\\
&\leq& n\|\bar p\|\big({\overline{\lambda} \delta t}\big)^{1/2}.
\end{eqnarray*}
One obtains the result by taking the $L^2$ expectation norm also on both sides of the inequalities involving $\epsilon_1$ and $\epsilon_2$.
\qquad\hspace{\stretch{1}}$ \Box$

Therefore, as the errors $\dE[|S_{j+1}-\wh{S}_{j+1}|^2]$ go to $0$ 
as the number of points in the discretization grids goes to infinity, we have the convergence of 
$\wt{P}(t)$ to ${P}(t)$ as long as the time grid step $\delta t$ also goes to $0$. 
Theorem \ref{th:Lip} also gives a convergence rate for $\|P(t)-\wt{P}(t)\|$, 
providing that $0\leq t< T\wedge T_{n+1}$. 
The convergence rate for the gain matrices is now straightforward from their definitions.
%
\begin{corollary}\label{cor:ErrK}
Under Assumption~\ref{hyp:infinite}, for all $0\leq t< T\wedge T_{n+1}$, one has
\begin{eqnarray*}
\lefteqn{\dE[\|K_{KB}(t)-\wt{K}(t)\|^2\ind{0\leq t\leq T\wedge T_{n+1}}]^{1/2}}\\
&\leq&\|C'(DD')^{-1}\|_\mathcal{S}\Big(\sum_{j=0}^{n-1}\bar\ell^{n-j}\bar\eta\dE[|S_{j+1}-\wh{S}_{j+1}|^2]^{1/2}\\
&&+\bar\eta \delta t+n\|\bar p\|{(\overline{\lambda}\delta t)}^{1/2}\Big).
\end{eqnarray*}
\end{corollary}
%
\subsection{Error derivation for the filtered trajectories}
We now turn to the estimation of the error between the exact KBF trajectory and our approximate one. We start with introducing some new notation. Let $b: \dR\times\dR^{2n_1}\rightarrow\dR^{2n_1}$ and $\wt{b}: \dR\times\dR^{2n_1}\rightarrow\dR^{2n_1}$ be defined by
\begin{eqnarray*}
b(t,z)&=&\left(
\begin{array}{cc}
A_{\theta(t)}&0\\
K_{KB}(t)C_{\theta(t)}&A_{\theta(t)}-K_{KB}(t)C_{\theta(t)}
\end{array}\right)z,\\
\wt{b}(t,z)&=&\left(
\begin{array}{cc}
A_{\theta(t)}&0\\
\wt K(t)C_{\theta(t)}&A_{\theta(t)}-\wt K(t)C_{\theta(t)}
\end{array}\right)z
\end{eqnarray*}
Let also $\sigma: \dR\rightarrow\dR^{2n_1\times(n_3+n_4)}$ and $\wt\sigma: \dR\rightarrow\dR^{2n_1\times(n_3+n_4)}$ be defined by
\begin{equation*}
\sigma(t)=\left(
\begin{array}{cc}
E_{\theta(t)}&0\\
0&K_{KB}(t)D_{\theta(t)}
\end{array}\right),
\end{equation*}
\begin{equation*}
\wt\sigma(t)=\left(
\begin{array}{cc}
E_{\theta(t)}&0\\
0&\wt K(t)D_{\theta(t)}
\end{array}\right).
\end{equation*}
Finally, set $W(t)=(w(t)',v(t)')'$, $X(t)=(x(t)',\wh{x}_{KB}(t)')'$ and $\wt{X}(t)=(x(t)',\wt{x}(t)')'$, so that the two processes $\{X(t), 0\leq t\leq T\}$ and $\{\wt{X}(t), 0\leq t\leq T\}$ have the following dynamics
\begin{equation*}
\left\{\begin{array}{l}
dX(t)=b(t,X_t)dt+\sigma(t)dW(t),\\
X(0)=(x(0)',\dE[x(0)]')',
\end{array}\right.
\end{equation*}
\begin{equation*}
\left\{\begin{array}{l}
d\wt{X}(t)=\wt{b}(t,\wt{X}_t)dt+\wt\sigma(t)dW(t),\\
\wt{X}(0)=(x(0)',\dE[x(0)]')'.
\end{array}\right.
\end{equation*}
The regularity properties of functions $b$, $\wt{b}$, $\sigma$ and $\wt{\sigma}$ are quite straightforward from their definition.
\begin{lemma}\label{lem:Lipbsig}
Under Assumption~\ref{hyp:infinite}, for all $0\leq t\leq T$ and $z,\wh{z}\in\dR^{2n_1}$, one has
\begin{eqnarray*}
|b(t,z)|\!\!&\!\!\leq\!\!&\!\!(\|A\|_\mathcal{S}+\|\bar p\|\|C\|_\mathcal{S}^2\|(DD')^{-1}\|_\mathcal{S})|z|,\\
|\wt{b}(t,z)|\!\!&\!\!\leq\!\!&\!\!(\|A\|_\mathcal{S}+\|\bar p\|\|C\|_\mathcal{S}^2\|(DD')^{-1}\|_\mathcal{S})|z|,\\
\|\sigma(t)\|_2\!\!&\!\!\leq\!\!&\!\!\|E\|_\mathcal{S}+\|\bar p\|\|C\|_\mathcal{S}\|(DD')^{-1}\|_2\|D\|_\mathcal{S},\\
\|\wt\sigma(t)\|_2\!\!&\!\!\leq\!\!&\!\!\|E\|_\mathcal{S}+\|\bar p\|\|C\|_\mathcal{S}\|(DD')^{-1}\|_\mathcal{S}\|D\|_\mathcal{S},\\
|b(t,z)-b(t,\wh{z})|\!\!&\!\!\leq\!\!&\!\!(\|A\|_\mathcal{S}+\|\bar p\|\|C\|_2^2\|(DD')^{-1}\|_\mathcal{S})|z-\wh{z}|,\\
|\wt{b}(t,z)-\wt{b}(t,\wh{z})|\!\!&\!\!\leq\!\!&\!\!(\|A\|_\mathcal{S}+\|\bar p\|\|C\|_\mathcal{S}^2\|(DD')^{-1}\|_\mathcal{S})|z-\wh{z}|,
\end{eqnarray*}
where $\bar p$ is the matrix defined in Lemma~\ref{lem0:Lip}.
\end{lemma}
\textit{Proof.} 
Upper bounds for $\|K_{KB}(t)\|_2$ and $\|\wt K(t)\|_2$ come from the upper bounds for $P_k(t)$ and $\wh{P}_k(t)$ given in Lemma~\ref{lem0:Lip}.
\qquad\hspace{\stretch{1}}$ \Box$

In particular, the processes $\{X(t), 0\leq t\leq T\}$ and $\{\wt{X}(t), 0\leq t\leq T\}$ are well defined and $\dE[\sup_{t\leq T}|X(t)|^2]$ and $\dE[\sup_{t\leq T}|\wt X(t)|^2]$ are finite, see e.g. \cite{KS91}. Set also $\Delta(t)=K_{KB}(t)-\wt K(t)$. In order to compare $X(t)$ and $\wt X(t)$, one needs first to be able to compare $b$ with $\wt{b}$ and $\sigma$ with $\wt{\sigma}$. The following result is straightforward from their definition.
\begin{lemma}\label{lem:Lipbsigt}
Under Assumption~\ref{hyp:infinite}, for all $0\leq t\leq T$ and $z\in\dR^{2n_1}$, one has
\begin{eqnarray*}
|b(t,z)-\wt b(t,z)|&\leq&2\|C\|_\mathcal{S}\|\Delta(t)\||z|,\\
\|\sigma(t)-\wt \sigma(t)\|_\mathcal{S}&\leq&\|D\|_\mathcal{S}\|\Delta(t)\|.
\end{eqnarray*}
\end{lemma}
We also need some bounds on the conditional moments of $\{X(t), 0\leq t\leq T\}$. Let $\{\mathcal{F}_t, 0\leq t\leq T\}$ be the filtration generated by the {semi-}Markov process $\{\theta(t), 0\leq t\leq T\}$, and $\dE_t[\cdot]=\dE[\cdot\ |\ \mathcal{F}_t]$. 
\begin{lemma}\label{lem:X4}
Under Assumption~\ref{hyp:infinite}, there exists a constant $c_2$ independent of the parameters of the system such that for $0\leq t \leq T$ one has
\begin{eqnarray*}
\lefteqn{\dE_T[\sup_{t\leq T\wedge T_{n+1}}|X(t)|^2]}\\
&\leq& 2c_2T(\|E\|_\mathcal{S}+\|\bar p\|\|C\|_\mathcal{S}\|(DD')^{-1}\|_\mathcal{S}\|D\|_\mathcal{S})^2\\
&&\times\exp(2T^2(\|A\|_\mathcal{S}+\|\bar p\|\|C\|_\mathcal{S}^2\|(DD')^{-1}\|_\mathcal{S})^2).
\end{eqnarray*}
\end{lemma}
\textit{Proof.} 
As $\{\theta(t), 0\leq t\leq T\}$ and the noise sequence $\{W(t), 0\leq t\leq T\}$ are independent, and the process $\{K_{KB}(t), 0\leq t\leq T\}$ is only dependent on $\{\theta(t), 0\leq t\leq T\}$ by construction, one has
\begin{eqnarray*}
\lefteqn{\dE_T[\sup_{u\leq t\wedge T\wedge T_{n+1}}|X(u)|^2]}\\
&\leq & 2\dE_T\Big[\sup_{u\leq t\wedge T\wedge T_{n+1}}\Big|\int_0^{u}\sigma(s)dW(s)\Big|^2\Big]\\
&&+2\dE_T\Big[\sup_{u\leq t\wedge T\wedge T_{n+1}}\Big|\int_0^{u}b(s,X(s))ds\Big|^2\Big]\\
&\leq&2c_2\dE_T\Big[\int_0^{T\wedge T_{n+1}}\big\|\sigma(s)\big\|^2ds\Big]\\
&&+2T\dE_T\Big[\int_0^{t\wedge T\wedge T_{n+1}}\big|b(s,X(s))\big|^2ds\Big],\\
\end{eqnarray*}
from convexity and Burkholder--Davis--Gundy inequalities, see e.g. \cite{KS91}, where $c_2$ is a constant independent of the parameters of the problem. From Lemma~\ref{lem:Lipbsig} one gets
\begin{eqnarray*}
\lefteqn{\dE_T[\sup_{u\leq t\wedge T\wedge T_{n+1}}|X(u)|^2]}\\
&\leq&2c_2T(\|E\|_\mathcal{S}+\|\bar p\|\|C\|_\mathcal{S}\|(DD')^{-1}\|_\mathcal{S}\|D\|_\mathcal{S})^2\\
&&+2T(\|A\|_\mathcal{S}+\|\bar p\|\|C\|_\mathcal{S}^2\|(DD')^{-1}\|_\mathcal{S})^2\\
&&\times\int_0^{t}{\dE_T[\sup_{u\leq s\wedge T\wedge T_{n+1}}|X(u)|^2]ds}.
\end{eqnarray*}
Finally, we use Gronwall's lemma to obtain
\begin{eqnarray*}
\lefteqn{\dE_T[\sup_{t\leq T\wedge T_{n+1}}|X(t)|^2]}\\
&\leq & 2c_2T(\|E\|_\mathcal{S}+\|\bar p\|\|C\|_\mathcal{S}\|(DD')^{-1}\|_\mathcal{S}\|D\|_\mathcal{S})^2\\
&&\times\exp(2T^2(\|A\|_\mathcal{S}+\|\bar p\|\|C\|_\mathcal{S}^2\|(DD')^{-1}\|_\mathcal{S})^2)
\end{eqnarray*}
which proves the result.
\qquad\hspace{\stretch{1}}$ \Box$

In the sequel, let $\overline{X}$ be the upper bound given by Lemma~\ref{lem:X4}:
\begin{eqnarray*}
\overline{X}&=&2c_2T(\|E\|_\mathcal{S}+\|\bar p\|\|C\|_\mathcal{S}\|(DD')^{-1}\|_\mathcal{S}\|D\|_\mathcal{S})^2\\
&&\times\exp(2T^2(\|A\|_\mathcal{S}+\|\bar p\|\|C\|_\mathcal{S}^2\|(DD')^{-1}\|_\mathcal{S})^2).
\end{eqnarray*}
We can now state and prove our convergence result.
\begin{theorem}\label{th:cv filter}
Under Assumption~\ref{hyp:infinite}, for $0\leq t \leq T$ one has
\begin{equation*}
\dE[|X(t)-\wt{X}(t)|^2\ind{0\leq t\leq T\wedge T_{n+1}}]\leq \overline{c}_1\exp(T\overline{c}_2),
\end{equation*}
with
\begin{eqnarray*}
\overline{c}_1&=&(2\|D\|_\mathcal{S}+8T\|C\|_\mathcal{S}^2\overline{X})\|C_{i}'(D_{i}D_{i}')^{-1}\|_\mathcal{S}\\
&&\times\Big(\sum_{j=0}^{n-1}\bar\ell^{n-j}\bar\eta\dE[|S_{j+1}-\wh{S}_{j+1}|^2]^{1/2}\\
&&+\bar\eta \delta t+n\|\bar p\|{(\overline{\lambda}\delta t)}^{1/2}\Big)^2,\\
\overline{c}_2&=&2T(\|A\|_\mathcal{S}+\|\bar p\|\|C\|_\mathcal{S}^2\|(DD')^{-1}\|_\mathcal{S})^2.
\end{eqnarray*}
\end{theorem}
\textit{Proof.} 
We follow the same lines as in the previous proof. As $\{\theta(t), 0\leq t\leq T\}$ and the noise sequence $\{W(t), 0\leq t\leq T\}$ are independent, and the processes $\{K_{KB}(t), 0\leq t\leq T\}$ and $\{\wt K(t), 0\leq t\leq T\}$ are only dependent on $\{\theta(t), 0\leq t\leq T\}$ by construction, one has
\begin{eqnarray*}
\lefteqn{\dE_T[|X(t)-\wt{X}(t)|^2\ind{0\leq t\leq T\wedge T_{n+1}}]}\\
&\leq & 2\dE_T\Big[\Big|\int_0^{t\wedge T\wedge T_{n+1}}\big(\sigma(s)-\wt\sigma(s)\big)dW(s)\Big|^2\Big]\\
&&+2\dE_T\Big[\Big|\int_0^{t\wedge T\wedge T_{n+1}}\big(b(s,X(s))-\wt b(s,\wt X(s))\big)ds\Big|^2\Big]\\
&\leq&2\dE_T\Big[\int_0^{t\wedge T\wedge T_{n+1}}\big\|\sigma(s)-\wt\sigma(s)\big\|^2ds\Big]\\
&&+2T\dE_T\Big[\int_0^{t\wedge T\wedge T_{n+1}}\big|b(s,X(s))-\wt b(s,\wt X(s))\big|^2ds\Big],\\
\end{eqnarray*}
from the isometry property of It\^o integrals and Cauchy--Schwartz inequality. 
From Lemmas~\ref{lem:Lipbsig}, \ref{lem:Lipbsigt} and Fubini one gets
\begin{eqnarray*}
\lefteqn{\dE_T[|X(t)-\wt{X}(t)|^2\ind{0\leq t\leq T\wedge T_{n+1}}]}\\
&\leq&2\|D\|_\mathcal{S}\int_0^{t\wedge T\wedge T_{n+1}}\big\|\Delta(s)\big\|^2ds\\
&&+2T\|C\|_\mathcal{S}^2\int_0^{t\wedge T\wedge T_{n+1}}\big\|\Delta(s)\big\|^2|\dE_T[|X(s)|^2]ds\\
&&+2T(\|A\|_\mathcal{S}+\|\bar p\|\|C\|_\mathcal{S}^2\|(DD')^{-1}\|_\mathcal{S})^2\\
&&\times\dE_T\Big[\int_0^{t\wedge T\wedge T_{n+1}}\big|X(s)-\wt X(s)\big|^2ds\Big]\\
&\leq&(2\|D\|_\mathcal{S}+8T\|C\|_\mathcal{S}^2\overline{X})\int_0^{t\wedge T\wedge T_{n+1}}\big\|\Delta(s)\big\|^2ds\\
&&+2T(\|A\|_\mathcal{S}+\|\bar p\|\|C\|_\mathcal{S}^2\|(DD')^{-1}\|_\mathcal{S})^2\\
&&\times\int_0^{t}\dE_T[\big|X(s)-\wt X(s)\big|^2\ind{0\leq s\leq T\wedge T_{n+1}}]ds\\
&\leq&\wt{c}_1+\wt{c}_2\int_0^{t}\dE_T[\big|X(s)-\wt X(s)\big|^2\ind{0\leq s\leq T\wedge T_{n+1}}]ds, 
\end{eqnarray*}
from Lemma~\ref{lem:X4}, with
\begin{eqnarray*}
\wt c_1&=&(2\|D\|_\mathcal{S}+8T\|C\|_\mathcal{S}^2\overline{X})\int_0^{t\wedge T\wedge T_{n+1}}\big\|\Delta(s)\big\|^2ds,\\
\wt c_2&=&2T(\|A\|_\mathcal{S}+\|\bar p\|\|C\|_\mathcal{S}^2\|(DD')^{-1}\|_\mathcal{S})^2.
\end{eqnarray*}
We use Gronwall's lemma to obtain
\begin{eqnarray*}
\dE_T[|X(t)-\wt{X}(t)|^2\ind{0\leq t\leq T\wedge T_{n+1}}]
&\leq &\wt c_1\exp(T\wt c_2),
\end{eqnarray*}
and conclude by taking the expectation on both sides and using Corollary~\ref{cor:ErrK} to bound $\dE[\wt c_1]$.
\qquad\hspace{\stretch{1}}$ \Box$

As a consequence of the previous result, $|\wh{x}_{KB}(t)-\wt{x}(t)|$ goes to $0$ almost surely as the number of points in the discretization grids goes to infinity.

\begin{rem}\label{rem-delayed-observation}
As noted in Remark \ref{rem-byproduct}, in the case of imperfect observation $\wt{S}_k$ of 
$S_k$, the errors $\dE[|\wt S_{j+1}-\wh{S}_{j+1}|^2]$ 
do not necessarily go to $0$ if $\theta$ is not instantaneously observed, 
however the errors are small when the time delays are small. 
The previous result implies that the filter performance deterioration 
is proportional to these errors. Acceptable performances can still be 
achieved in applications where $\theta$ is not instantaneously observed.
\end{rem}

%
\section{Numerical example}
\label{sec-example}
We now illustrate our results on a magnetic suspension system presented in \cite{ECosta99ieee}. 
The system is a laboratory device that consists of a coil whose voltage is 
controlled by a rather simple (non-reliable) pulse-width modulation system, 
and sensors for position of a suspended metallic sphere and for the 
coil current. The model around the origin without jumps and noise is
in the form $\dot x(t)=A x(t)+B u(t)$, $y(t)=C x(t)$, with
\begin{equation*}%
A=\left(\begin{array}{ccc}
0&1&0\\
1750&0&-34.1\\
0& 0 & -0.0383
\end{array}\right),\qquad
B=\left(\begin{array}{c}
0\\
0\\
1.9231
\end{array}\right),
\end{equation*}
\begin{equation*}
C=\left(\begin{array}{ccc}
1&0&0\\
0&0&1
\end{array}\right).
\end{equation*}
The components of vector $x(t)$ are the position of the sphere, its
speed and the coil current. The coil voltage $u(t)$  is controlled
using a stabilizing state feedback control,
leading to the closed loop dynamics $\dot x(t)=A_1 x(t)$,
\begin{equation*}%
A_1=\left(\begin{array}{ccc}
0&1&0\\
1750&0&-34.1\\
4360.2&104.2&-84.3
\end{array}\right).
\end{equation*}
We consider the realistic scenario where the system may be operating in
normal mode $\theta=1$ or in critical failure $\theta=2$
due e.g. to faults in the pulse-width modulation system, which is included in the model
by making $B_2=0$, leading to the closed loop dynamics $\dot x(t)=A_2 x(t)$
with $A_2=A$.
Although it is natural is to consider that the system starts in normal
mode a.s. and never recovers from a failure,  
we want to compare the performance of the proposed filter with
the LMMSE \cite{FC10} that requires {a true Markov chain with} positive probabilities
for all modes at all times, then we relax the problem by setting
the initial distribution $\pi(0)=(0.999, 0.001)$ and {the transition rates matrix}
\begin{equation*}
\Lambda=\left(\begin{array}{cc}
-20&20\\
0.1&-0.1
\end{array}\right)
\end{equation*}
with the interpretation that the recovery from failure mode is
relatively slow. 

In the overall model  Eq.~\eqref{mjls} we set $C_1=C_2=C$ and we also consider that $x(0)$ is
normally distributed with mean
 $\mathbb{E}[x(0)]=(0.001,0,0)'$ and variance $Var(x(0))=I_3$,
\begin{equation*}
E_1\!=\!E_2\!=\!\left(\begin{array}{ccc}
1&0.2&-1.9\\
-0.1&1.4&-0.3\\
0.1&0.5&1
\end{array}\right)\!\!,\ 
D_1\!=\!D_2\!=\!\left(\begin{array}{cc}
1&0\\
0&1
\end{array}\right)\!\!,
\end{equation*}
so that only the position of the sphere and the coil current are
measured through some noise. Speed is not observed. 
It is worth mentioning that the system is not mean square 
stable, so that the time horizon $T$ is usually short for the 
trajectory to stay close to the origin and keep the linearized model valid; 
we can slightly increase the horizons during simulations for academic purposes only.
%
\subsection{Markovian linear minimum mean squares estimator}
Fragoso and Costa proposed in \cite{FC10} the so-called Markovian linear minimum
mean squares estimator (LMMSE) for MJLS with finite state space Markov chains. Under Assumption~\ref{hyp:finite},
 the equation of the filter is
\begin{eqnarray*}
d\hat{x}_{FC}(t)&=&A_{\theta(t)}\hat{x}_{FC}(t)dt\\
&&+K_{FC}(\theta(t),t)(dy(t)-C_{\theta(t)}\hat{x}_{FC}(t)dt),
\end{eqnarray*}
for $0\leq t\leq T$, with initial condition $\hat{x}_{FC}(0)=\mathbb{E}[x(0)]$ and gain matrices
\begin{equation*}
K_{FC}(i,t)=P_{FC}(i,t)C_{i}'(D_{i}D_{i}'\pi_i(t))^{-1},
\end{equation*}
where $\pi_i(t)=\mathbb{P}(\theta(t)=i)=(\pi(0)\exp(t\Lambda))_i$ and
$\{P_{FC}(i,t), 0\leq t\leq T\}$ satisfies the system of matrix differential equation
\begin{eqnarray*}
dP_{FC}(i,t)&=&\big(A_{i}P_{FC}(i,t)+P_{FC}(i,t)A_{i}'\\
&&+\sum_{j=1}^N  P_{FC}(j,t)\Lambda_{ji}+E_{i}E_{i}'\pi_i(t)\\
&&-P_{FC}(i,t)C_{\theta(t)}'(D_{\theta(t)}D_{\theta(t)}'\pi_i(t))^{-1}\\
&&\times C_{\theta(t)}P_{FC}(i,t)\big)dt,\\
P_{FC}(i,0)&=&Var(x(0))\pi_i(0).
\end{eqnarray*}
The matrices $\{P_{FC}(i,t), 0\leq t\leq T, i\in\mathcal{S}\}$ and $\{K_{FC}(i,t), 0\leq t\leq T, i\in\mathcal{S}\}$ depend only on
the law of $\{\theta(t), 0\leq t\leq T\}$ and not on its current value.
Therefore they can be computed off line on a discrete time grid and
stored but it is sub-optimal compared to the KBF.
%
\subsection{Approximate filter by quantization}
We start with the quantized discretization of the inter-jump times $\{S_n\}$ of the Markov chain $\{\theta(t), 0\leq t\leq T\}$. We use the CLVQ algorithm described for instance in \cite{pages98}. Table~\ref{tab:quant error JT} gives the error $\dE[|S_1-\wh{S}_1|^2\ |\ \theta(0)=i]^{1/2}$ for $i=1,2$ computed with $10^6$ Monte Carlo simulations for an increasing number of discretization points. This illustrates the convergence of Theorem~\ref{th:quantize}: the error decreases as the number of points increases. The variance of the first jump time in mode $2$ is much higher than in mode $1$ which accounts for the different scales in the errors.
\begin{figure}[t]
\centerline{\includegraphics[height=5cm]{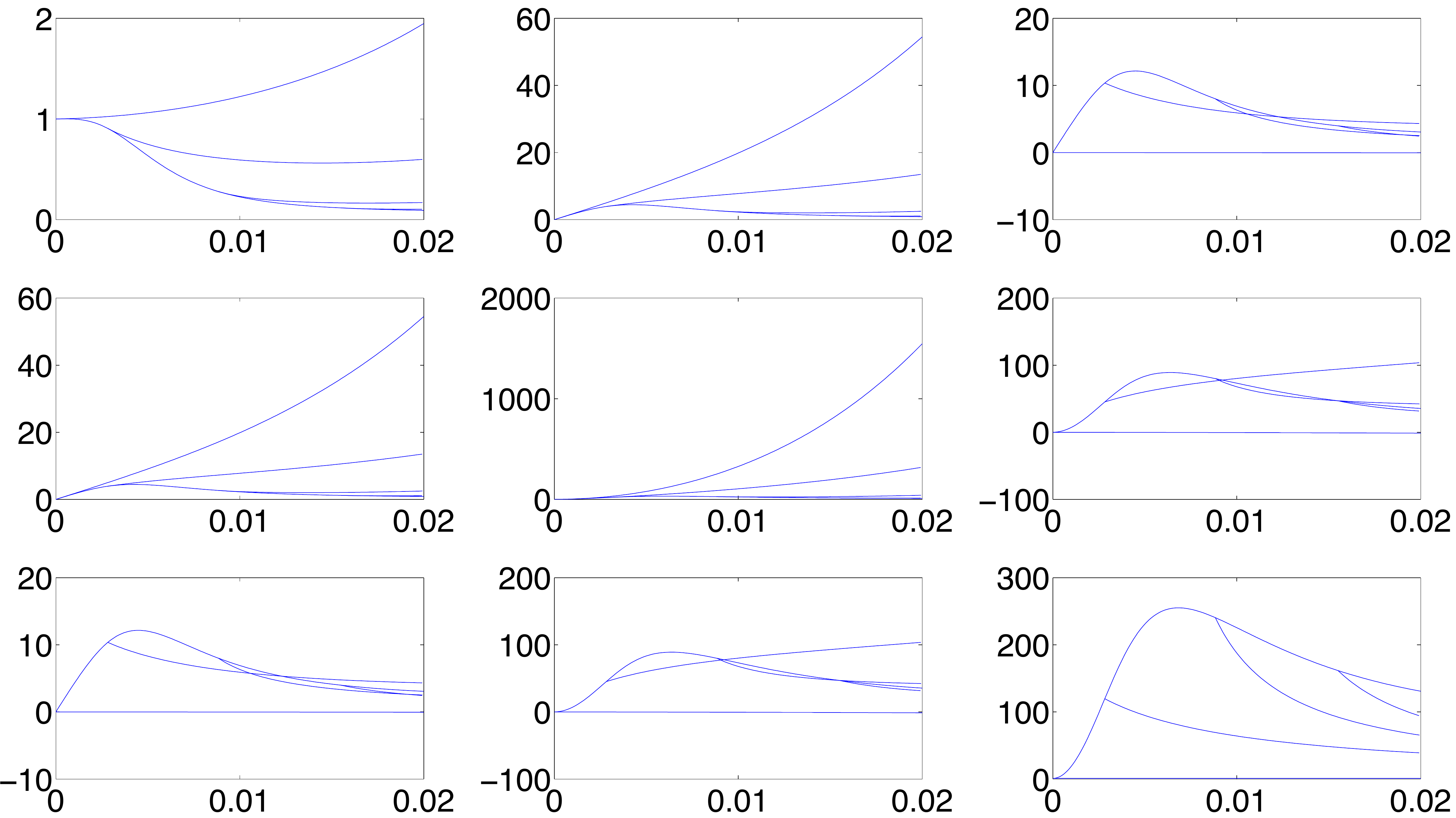}}
\caption{Pre-computed tree of solutions with $10$ grid points.}
\label{tree10}
\end{figure}
\begin{figure}[t]
\centerline{\includegraphics[height=5cm]{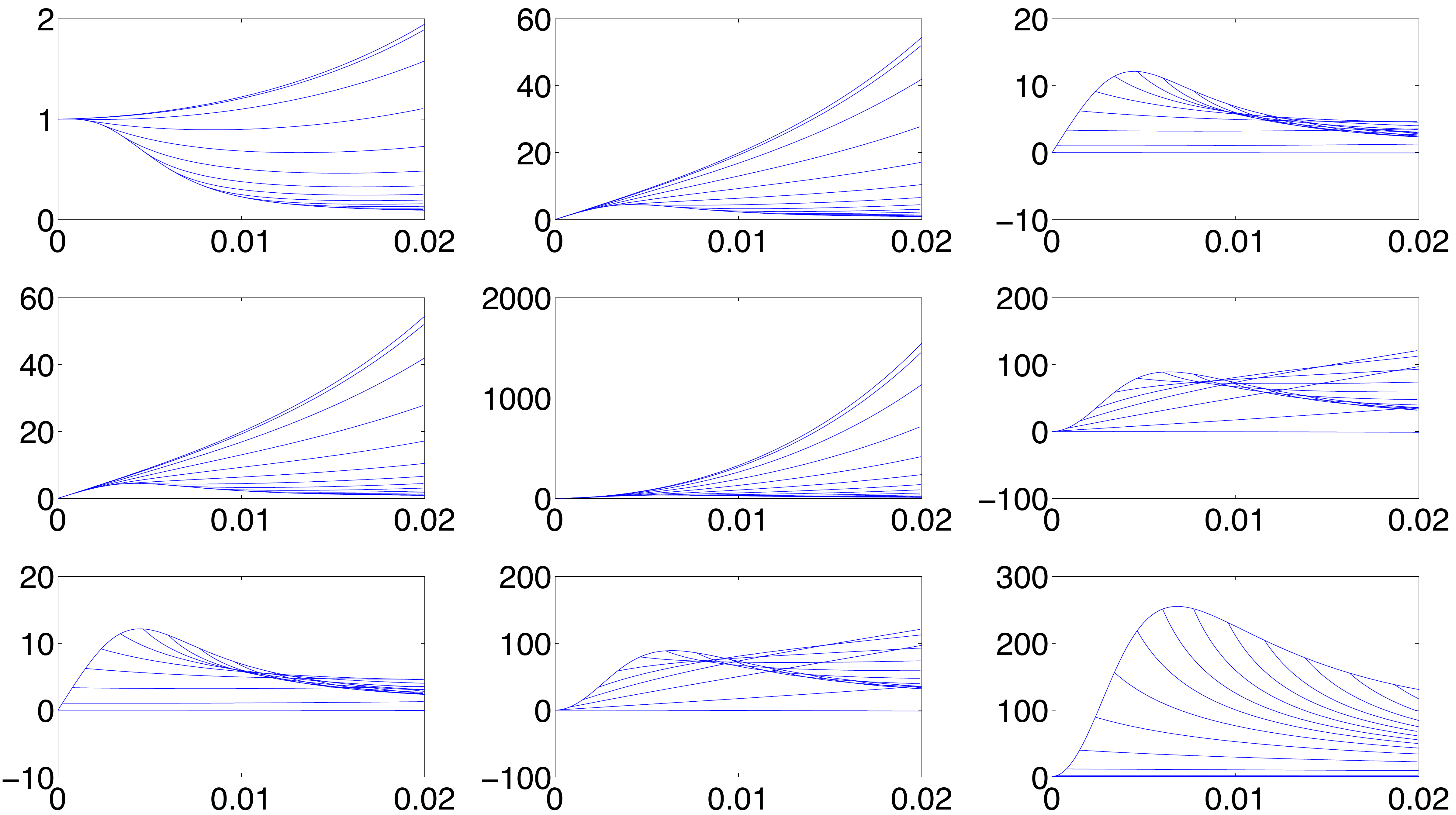}}
\caption{Pre-computed tree of solutions with $50$ grid points.}
\label{tree50}
\end{figure}
\begin{table}[t]
\footnotesize
\begin{center}
\begin{tabular}{ccc}
\hline
Number of grid points&Error for $\theta(0)=1$&Error for $\theta(0)=2$\\
\hline
10& 5.441$\times 10^{-3}$&1017$\times 10^{-3}$\\
50&1.585$\times 10^{-3}$&357.5$\times 10^{-3}$\\
100&0.753$\times 10^{-3}$&175.2$\times 10^{-3}$\\
500&0.173$\times 10^{-3}$&36.22$\times 10^{-3}$\\
1000&0.100$\times 10^{-3}$& 23.35$\times 10^{-3}$\\
\hline
\end{tabular}
\caption{Quantization error for the first jump time depending on the number of points in the discretization grid and the value of the starting point of the Markov chain.}
\label{tab:quant error JT}
\end{center}
\end{table}

The second step consists in solving the Riccati equation (\ref{eq Ric}) for all possible trajectories of $\{\theta(t), 0\leq t\leq T\}$ with inter-jump times in the quantization grids and up to the computation horizon $T=0.02$. Namely, we compute the trajectories \{$\wh{P}_k(t), 0\leq t\leq T\}$. We chose a regular time grid with time step $\delta t=10^{-4}$. For technical reasons related to the selection of branches, the time horizon $T$ is added in each grid. One thus obtains a tree of pre-computed branches that are solutions of Eq. (\ref{eq Ric}), the branching times being the quantized jump times.
Figures~\ref{tree10} and \ref{tree50} show the pre-computed trees of solutions component-wise for $10$ and $50$ points respectively in the quantization grids. Note the very different scales of the coordinates.
The number of grid points that are actually used (quantized points below the horizon $T$) are given in Table~\ref{tab:points below horizon}  for each original quantization grid size, together with the resulting number of pre-computed branches. 
\begin{table}[t]
\footnotesize
\begin{center}
\begin{tabular}{cccr}
\hline
Number of &Points below &Points below&Number of\\
grid points&horizon&horizon&branches\\
&for $\theta(0)=1$& for $\theta(0)=2$&\\
\hline
10& 4&1&7\\
50&14&1&17\\
100&33&1&36\\
500&161&2&7763\\
1000&319&3&603784\\
\hline
\end{tabular}
\caption{Number of grid points actually used and corresponding number of pre-computed branches depending on the initial number of points in the discretization grid.}
\label{tab:points below horizon}
\end{center}
\end{table}
The number of pre-computed branches grows exponentially fast when we take into account more grid points. Time taken to pre-compute the branches grows accordingly. In this example, the number of points used in mode $2$ is low, therefore the number of branches remains tractable.

\begin{figure}[t]
\centerline{\includegraphics[height=4cm]{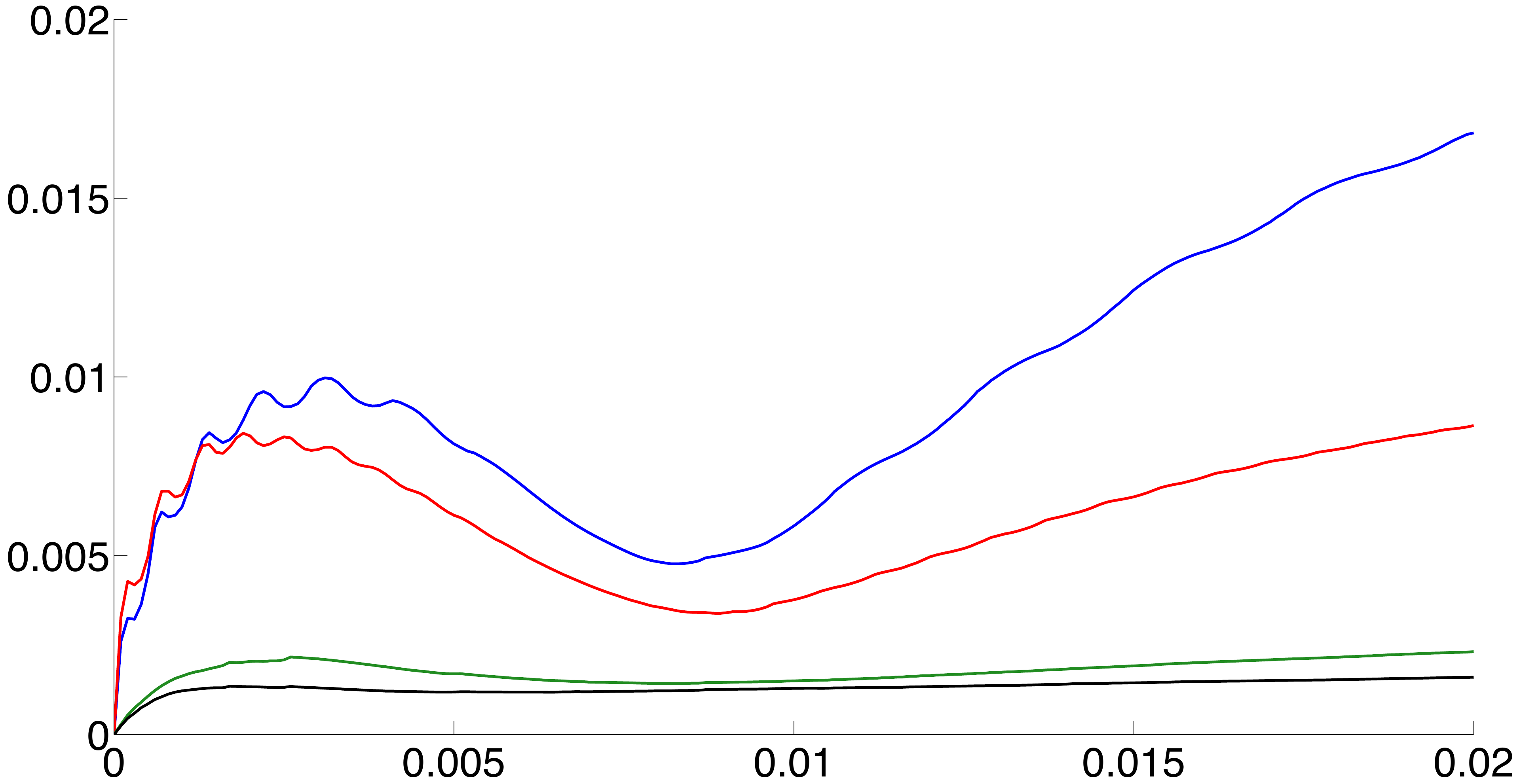}}
\caption{Average relative error between the solution of Riccati equation and its approximation, 
from top to bottom: blue: 50 points, red: 100 points, green: 500 point, black: 1000 points in the quantization grids.}
\label{ric-error}
\end{figure}

To compute the filtered trajectory in real time, one starts with the approximation of the solution of Eq. (\ref{eq Ric}). The first branch corresponds to the pre-computed branch starting at time $0$ from $\theta(0)$. When the first jump occurs, one selects the nearest neighbor of the jump time in the quantization grid and the corresponding pre-computed branch, and so on for the following jumps.
Figure \ref{ric-error} shows the mean of the relative error between the solution of Eq (\ref{eq Ric}) and its approximation (for the matrix norm 2) for given numbers of points in the quantization grids and $10^5$ Monte Carlo simulations. Again, it illustrates how the accuracy of the approximation increases with the number of points in the quantization grids.

\begin{figure}[t]
\centerline{\includegraphics[height=4cm]{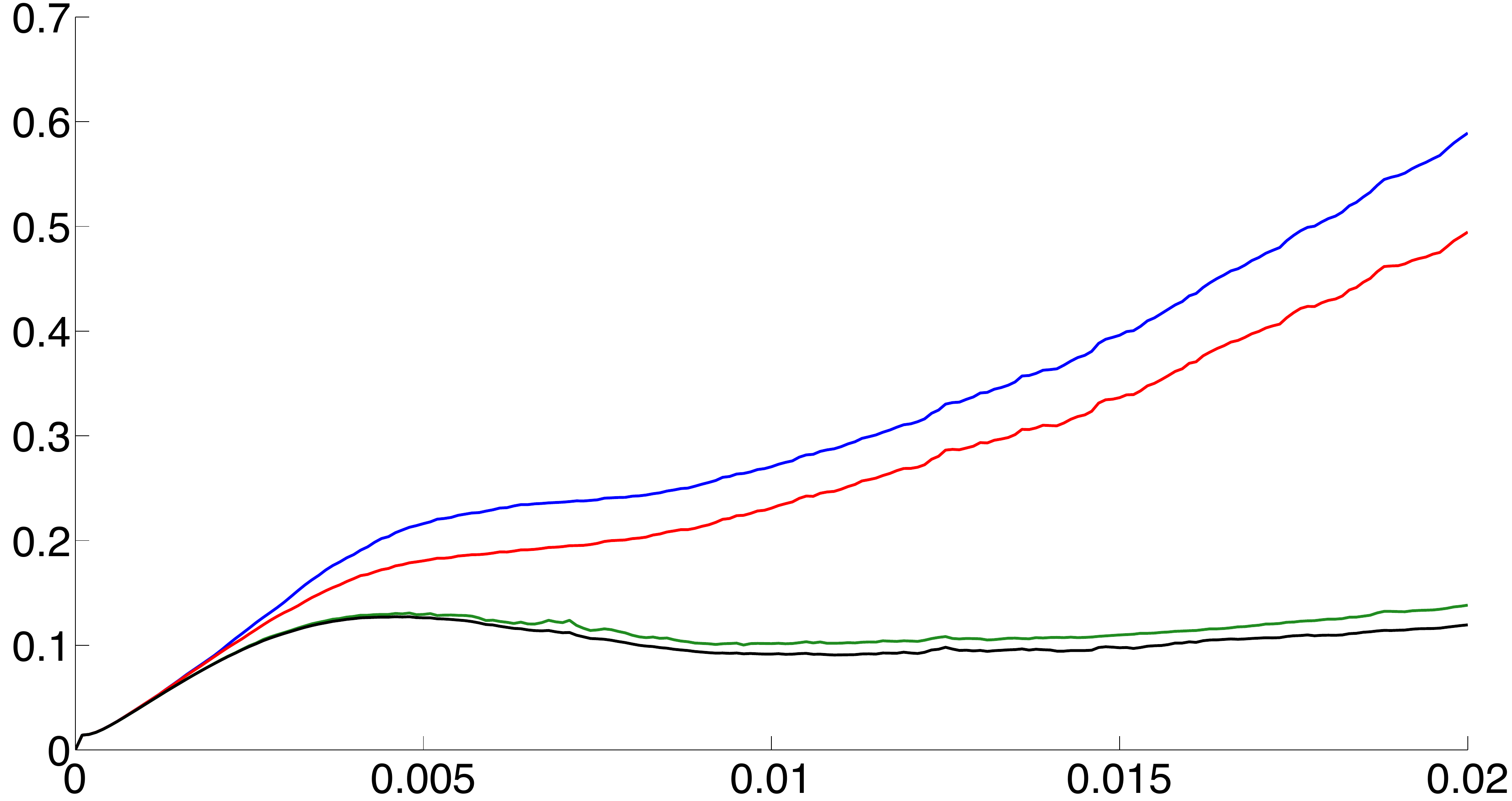}}
\caption{$L^2$ norm of the  difference between $\wh{x}_{KB}$ and its quantized approximation $\wt{x}$, from top to bottom: blue: 50 points, red: 100 points, green: 500 point, black: 1000 points in the quantization grids.}
\label{filter-error}
\end{figure}

Finally, the real-time approximation of Eq (\ref{eq Ric}) is plugged into the filtering equations to obtain an approximate KBF.
Figure \ref{filter-error} shows the mean $L^2$ distance between the real KBF $\{\wh{x}_{KB}(t), 0\leq t\leq T\}$ and its approximation $\{\wt{x}, 0\leq t\leq T\}$ following our procedure for an increasing number of points in the quantization grids and for $10^5$ Monte Carlo simulations.
%
\subsection{Comparison of the filters}
For each filter, we ran $10^5$ Monte Carlo simulations and computed the mean of the following error between the real trajectory $\{x(t), 0\leq t\leq T\}$ and the filtered trajectory $\{\hat{x}(t), 0\leq t\leq T\}$ for all of the three filters presented above, the exact Kalman--Bucy filter being the reference.
\begin{equation*}
\int_0^T\Big(\big(x_1(t)-\hat{x}_1(t)\big)^2+\big(x_2(t)-\hat{x}_2(t)\big)^2+\big(x_3(t)-\hat{x}_3(t)\big)^2\Big) dt.
\end{equation*}
Table~\ref{tab:filter_error} gives this error for given numbers of points in the quantization grids. Of course only the error for the approximate filter changes with the quantization grids. Note that our approximate filter is very close to the KBF and performs better than the LMMSE for as little as $10$ points in the quantization grids corresponding to $7$ precomputed branches.
\begin{table}[t]
\footnotesize
\begin{center}
\begin{tabular}{cccc}
\hline
Number of grid points&Error for &Error for &Error for \\
&KBF&approximate filter&LMMSE\\
\hline
10&3.9244&3.9634&3.9850\\
50&3.9244&3.9254&3.9850\\
100&3.9244&3.9246&3.9850\\
500&3.9244&3.9244&3.9850\\
1000&3.9244&3.9244&3.9850\\
\hline
\end{tabular}
\caption{Average error for the different filters depending on the number of points in the quantization grids, 
considering horizon $T=0.02$.}
\label{tab:filter_error}
\end{center}
\end{table}
We also ran our simulations with longer horizons. The performance of the filters is given 
in Table~\ref{tab:filter_errorLong} and illustrate that our filter can still perform good with a longer horizon. 
Note that the computations of the LMMSE is impossible from an horizon of $0.4$ on because the estimated 
state space reaches too high values very fast, and they are treated as infinity numerically. 
From an horizon of $0.8$ on, all computations are impossible because the system is not mean square 
stable, as we explained before.
\begin{table}[t]
\footnotesize
\begin{center}
\begin{tabular}{crrrrr}
\hline
$T$&Grid &Branches&Error for &Error for &Error for \\
& points&&KBF&approx. filter&LMMSE\\
\hline
0.1&10&12&376.3&425.6&812.5\\
0.1&50&110&376.3&379.1&812.5\\
0.1&100&3519&376.3&376.6&812.5\\
\hline
0.2&10&14&8597&10610&13260\\
0.2&50&2832&8597&9715&13260\\
\hline
0.3&10&14&2.325$\times 10^4$&4.893$\times 10^6$&3.023$\times 10^5$\\
0.3&50&11248&2.325$\times 10^4$&4.141$\times 10^6$&3.023$\times 10^5$\\
\hline
0.4&10&14&4.913$\times 10^4$&4.663$\times 10^{10}$&NaN\\
0.4&50&50049&4.913$\times 10^4$&2.102$\times 10^{10}$&NaN\\
\hline
\end{tabular}
\caption{Average error for the different filters depending on the horizon, the number of points in the quantization grids and the number of branches.}
\label{tab:filter_errorLong}
\end{center}
\end{table}
%
\section{Conclusion}
\label{sec-conclusion}
We have presented a filter for state estimation of {s}MJLS relying on
discretization
by quantization of the {semi-}Markov chain and solving a finite number of filtering
Riccati equations.
The difference between the approximated Riccati solution $\widetilde P(t)$
and the actual Riccati solution $P(t)$ has been studied
and we have shown in  Theorem \ref{th:Lip} that it converges to zero in average
when the number of points in the discretization grid goes to infinity;
a convergence rate is also provided, allowing to find a convergence
rate for the gain matrices, see Corollary \ref{cor:ErrK}.
Based on this result, and on an upper bound for the conditional second moment of the KBF
that is derived in Lemma \ref{lem:X4}, we have obtained the main
convergence result in Theorem \ref{th:cv filter}, which implies
convergence to zero of $\mathbb{E}|x_{KB}(t)-\widetilde x(t)|^2$, so that
$\widetilde x(t)$ approaches $x_{KB}(t)$ almost surely as the number of
grid points goes to infinity. 
Applications in which $\theta$ is not instantaneously observed can 
also benefit from the proposed filter, however it may not completely recover the 
performance of the  KBF as explained in Remark \ref{rem-delayed-observation}.
The algorithm has been applied to a real-world system and performed
almost as well as the KBF with a small grid of $10$ points.

Although the proposed filter can be pre-computed, the number of
branches of the Riccati equation grows exponentially with the time horizon $T$,
making the pre-computation time too high in some cases. One exception comprises
systems with no more than one fast mode (high transition rates),
because in such a situation the slow modes do not branch much and the
number of branches grows in an almost linear fashion with $T$ as long as 
the probability of the slow mode to jump before $T$ remains small.
Examples of applications coping with this setup,
which can benefit from the proposed filter, are systems with
small probability of failure and quick recovery (the failure mode is fast),
or a variable number of permanent failures (the normal mode is fast),
with web-based control as a fertile field of applications.
For general systems, one possible way out of this cardinality issue
is to use a rolling-horizon scheme where the approximate gains are
pre-computed in small batches during the system operation and sent
to the controller memory. Another approach could be to quantize directly the sequence $\{S_k, P_k(S_k)\}$ thus keeping the number of branches at a fixed number, allowing for general transition rate matrices and longer horizons in terms of the number of jumps. However this approach suffers from a curse of dimensionality as the quantization error goes to zero with slower and slower rate as the dimension of the process goes higher, see Theorem~\ref{th:quantize}.

Future work will look into a rolling-horizon implementation scheme,
implementation issues and different compositions of the
KBF/LMMSE, for instance using time-delayed solutions of the KBF that can be computed during
the system operation as a measure for discarding unnecessary branches.
Alternative schemes for discretization/quantization and selection of
the appropriate
pre-computed solutions can be pursued, seeking to reduce the computational load
of the current algorithm while preserving the quality of the estimate.

\section*{Acknowledgment}
This work was supported by FAPESP Grant 13/19380-8, 
CNPq Grants 306466/2010 and 311290/2013-2, USP-COFECUB  
 FAPESP/FAPs/INRIA/INS2i-CNRS Grant 13/50759-3, Inria Associate team CDSS and ANR Grant Piece ANR-12-JS01-0006.

\bibliographystyle{acm}
\bibliography{biblio-ric}

\begin{thebibliography}{10}

\bibitem{Anderson79}
{\sc Anderson, B. D.~O., and Moore, J.~B.}
\newblock {\em Optimal Filtering}, first~ed.
\newblock Prentice-Hall, London, 1979.

\bibitem{bally03}
{\sc Bally, V., and Pag{\`e}s, G.}
\newblock A quantization algorithm for solving multi-dimensional discrete-time
  optimal stopping problems.
\newblock {\em Bernoulli 9}, 6 (2003), 1003--1049.

\bibitem{bally05}
{\sc Bally, V., Pag{\`e}s, G., and Printems, J.}
\newblock A quantization tree method for pricing and hedging multidimensional
  {A}merican options.
\newblock {\em Math. Finance 15}, 1 (2005), 119--168.

\bibitem{brandejsky12}
{\sc {B}randejsky, A., de~{S}aporta, B., and {D}ufour, F.}
\newblock Numerical methods for the exit time of a piecewise-deterministic
  {M}arkov process.
\newblock {\em Adv. in Appl. Probab. 44}, 1 (2012), 196--225.

\bibitem{brandejsky13}
{\sc Brandejsky, A., de~{S}aporta, B., and Dufour, F.}
\newblock Optimal stopping for partially observed piecewise-deterministic
  {M}arkov processes.
\newblock {\em Stochastic Process. Appl. 123}, 8 (2013), 3201--3238.

\bibitem{Campo91}
{\sc Campo, L., Mookerjee, P., and Bar-Shalom, Y.}
\newblock State estimation for systems with sojourn-time-dependent {M}arkov
  model switching.
\newblock {\em Automatic Control, IEEE Transactions on 36}, 2 (Feb 1991),
  238--243.

\bibitem{ECosta99ieee}
{\sc Costa, E.~F., Oliveira, V.~A., and Vargas, J.~B.}
\newblock Digital implementation of a magnetic suspension control system for
  laboratory experiments.
\newblock {\em IEEE Transactions on Education 42\/} (1999), 315 -- 322.

\bibitem{Costa11autom_filter}
{\sc Costa, O.~L., and Benites, G.~R.}
\newblock Linear minimum mean square filter for discrete-time linear systems
  with {Markov} jumps and multiplicative noises.
\newblock {\em Automatica 47}, 3 (2011), 466 -- 476.

\bibitem{CostaFragosoMarques05}
{\sc Costa, O. L.~V., Fragoso, M.~D., and Marques, R.~P.}
\newblock {\em Discrete-Time Markovian Jump Linear Systems}.
\newblock Springer-Verlag, New York, 2005.

\bibitem{CostaFragosoTodorov}
{\sc Costa, O. L.~V., Fragoso, M.~D., and Todorov, M.~G.}
\newblock {\em Continuous-Time {Markov} Jump Linear Systems}.
\newblock Springer, Berlin, Heidelberg, 2013.

\bibitem{saporta10}
{\sc de~Saporta, B., Dufour, F., and Gonzalez, K.}
\newblock Numerical method for optimal stopping of piecewise deterministic
  {M}arkov processes.
\newblock {\em Ann. Appl. Probab. 20}, 5 (2010), 1607--1637.

\bibitem{saporta12}
{\sc de~{S}aporta, B., {D}ufour, F., {Z}hang, H., and {E}legbede, C.}
\newblock {O}ptimal stopping for the predictive maintenance of a structure
  subject to corrosion.
\newblock {\em Proceedings of the {I}nstitution of {M}echanical {E}ngineers,
  {P}art {O}: {J}ournal of {R}isk and {R}eliability 226}, 2 (2012), 169--181.

\bibitem{doVal99JE}
{\sc do~Val, J., and Basar, T.}
\newblock Receding horizon control of jump linear systems and a macroeconomic
  policy problem.
\newblock {\em Journal of Economic Dynamics \& Control 23\/} (1999),
  1099--1131.

\bibitem{Dragan2009}
{\sc Dragan, V., Morozan, T., and Stoica, A.~M.}
\newblock {\em Mathematical methods in robust control of discrete-time linear
  stochastic systems}.
\newblock Springer, 2009.

\bibitem{Dragan2013}
{\sc Dragan, V., Morozan, T., and Stoica, A.~M.}
\newblock {\em Mathematical Methods in Robust Control of Linear Stochastic
  Systems}, 2nd~ed.
\newblock Springer, 2013.

\bibitem{FC10}
{\sc Fragoso, M., and Costa, O. L.~V.}
\newblock A separation principle for the continuous-time {LQ}-problem with
  {M}arkovian jump parameters.
\newblock {\em IEEE Transactions on Automatic Control 55}, 12 (2010),
  2692--2707.

\bibitem{Geromel09}
{\sc Geromel, J., Gonçalves, A., and Fioravanti, A.}
\newblock Dynamic output feedback control of discrete-time {M}arkov jump linear
  systems through linear matrix inequalities.
\newblock {\em SIAM Journal on Control and Optimization 48}, 2 (2009),
  573--593.

\bibitem{gray98}
{\sc Gray, R.~M., and Neuhoff, D.~L.}
\newblock Quantization.
\newblock {\em IEEE Trans. Inform. Theory 44}, 6 (1998), 2325--2383.
\newblock Information theory: 1948--1998.

\bibitem{Hou06}
{\sc Hou, Z., Luo, J., Shi, P., and Nguang, S.~K.}
\newblock Stochastic stability of ito differential equations with
  semi-{M}arkovian jump parameters.
\newblock {\em Automatic Control, IEEE Transactions on 51}, 8 (Aug 2006),
  1383--1387.

\bibitem{Huang13}
{\sc Huang, J.}
\newblock {\em Analysis and Synthesis of Semi-{M}arkov Jump Linear Systems and
  Networked Dynamic Systems}.
\newblock PhD thesis, University of Victoria, 2013.

\bibitem{Huang13IJRNC}
{\sc Huang, J., and Shi, Y.}
\newblock Stochastic stability and robust stabilization of semi-{M}arkov jump
  linear systems.
\newblock {\em International Journal of Robust and Nonlinear Control 23}, 18
  (2013), 2028--2043.

\bibitem{Jazwinski70}
{\sc Jazwinski, A.~H.}
\newblock {\em Stochastic Processes and Filtering Theory}.
\newblock Academic Press, 1970.

\bibitem{Kalman60}
{\sc Kalman, R.}
\newblock A new approach to linear ltering and prediction problems.
\newblock {\em J. Basic Engineering 82}, 1 (1960), 35--45.

\bibitem{Kalman61}
{\sc Kalman, R., and Bucy, R.}
\newblock New results in linear ltering and prediction theory.
\newblock {\em J. Basic Engineering 83\/} (1961), 95--108.

\bibitem{KS91}
{\sc Karatzas, I., and Shreve, S.~E.}
\newblock {\em Brownian motion and stochastic calculus}, second~ed., vol.~113
  of {\em Graduate Texts in Mathematics}.
\newblock Springer-Verlag, New York, 1991.

\bibitem{Kumar86}
{\sc Kumar, P.~R., and Varaiya, P.}
\newblock {\em Stochastic Systems: Estimation, Identification, and Adaptive
  Control}.
\newblock Prentice-Hall, 1986.

\bibitem{pages98}
{\sc Pag{\`e}s, G.}
\newblock A space quantization method for numerical integration.
\newblock {\em J. Comput. Appl. Math. 89}, 1 (1998), 1--38.

\bibitem{pages05}
{\sc Pag{\`e}s, G., and Pham, H.}
\newblock Optimal quantization methods for nonlinear filtering with
  discrete-time observations.
\newblock {\em Bernoulli 11}, 5 (2005), 893--932.

\bibitem{pages04b}
{\sc Pag{\`e}s, G., Pham, H., and Printems, J.}
\newblock An optimal {M}arkovian quantization algorithm for multi-dimensional
  stochastic control problems.
\newblock {\em Stoch. Dyn. 4}, 4 (2004), 501--545.

\bibitem{pages04}
{\sc Pag{\`e}s, G., Pham, H., and Printems, J.}
\newblock Optimal quantization methods and applications to numerical problems
  in finance.
\newblock In {\em Handbook of computational and numerical methods in finance}.
  Birkh\"auser Boston, Boston, MA, 2004, pp.~253--297.

\bibitem{Schwartz03}
{\sc Schwartz, C.}
\newblock {\em Control of semi-{M}arkov jump linear systems with application to
  the bunch-train cavity interaction}.
\newblock PhD thesis, Northwestern University, 2003.

\bibitem{Siqueira04}
{\sc Siqueira, A. A.~G., and Terra, M.~H.}
\newblock Nonlinear and markovian ${H}_{\infty}$-controls of underactuated
  manipulators.
\newblock {\em IEEE Transactions on Control System Technology 12\/} (2004),
  811--826.

\bibitem{Sworder83}
{\sc Sworder, D.~D., and Rogers, R.~O.}
\newblock An {LQ}-solution to a control problem associated with a solar thermal
  central receiver.
\newblock {\em IEEE Transactions on Automatic Control 28}, 10 (1983), 971--978.

\end{thebibliography}

\end{document}